\documentclass[a4paper,oneside,british,12pt]{amsart}
\usepackage[T1]{fontenc}
\usepackage{textcomp}
\usepackage[latin9]{inputenc}
\usepackage{url}
\usepackage{xcolor}
\usepackage{mathtools}
\usepackage{enumitem}
\usepackage{multicol}
\usepackage{amsbsy}
\usepackage{amstext}
\usepackage{amsthm}
\usepackage{amssymb}
\usepackage{graphicx}
\usepackage{rotating}

\makeatletter

\providecommand{\tabularnewline}{\\}

\numberwithin{equation}{section}
\numberwithin{figure}{section}
\theoremstyle{remark}
\newtheorem*{rem*}{\protect\remarkname}
\newlist{casenv}{enumerate}{4}
\setlist[casenv]{leftmargin=*,align=left,widest={iiii}}
\setlist[casenv,1]{label={{\itshape\ \casename} \arabic*.},ref=\arabic*}
\setlist[casenv,2]{label={{\itshape\ \casename} \roman*.},ref=\roman*}
\setlist[casenv,3]{label={{\itshape\ \casename\ \alph*.}},ref=\alph*}
\setlist[casenv,4]{label={{\itshape\ \casename} \arabic*.},ref=\arabic*}

\@ifundefined{date}{}{\date{}}
\makeatother

\newcommand*{\ran}{\operatorname{ran}}
\newcommand*{\dom}{\operatorname{dom}}

\usepackage{listings}
\usepackage{babel}
\addto\captionsbritish{\renewcommand{\casename}{Case}}
\addto\captionsbritish{\renewcommand{\remarkname}{Remark}}
\addto\captionsenglish{\renewcommand{\casename}{Case}}
\addto\captionsenglish{\renewcommand{\remarkname}{Remark}}
\providecommand{\casename}{Case}
\providecommand{\remarkname}{Remark}

\begin{document}
\date{}
\title[On Lewis Carroll's Work on Logic]{On the Algebraic Structure of Lewis Carroll's Logic\\
{\footnotesize\emph{Behind the Looking Glass: A Mathematical Perspective}}}
\author{Sebastian Franz \& Rainer Picard}

\maketitle
\tableofcontents{}

\setcounter{section}{-1}

\section{Introduction}

\subsection{\protect The idea}\label{subsec:The-idea}

Charles Ludwidge Dodgeson, better known by his literary pen name `Lewis
Carroll', has a number of well known and cherished literary works
to his credit, for example ``Alice in Wonderland'' and ``Through
the Looking Glass''. Apart from this he is known to have a strong
professional interest in mathematics and logic with an inclination
towards an educational and recreational perspective. As it turns out,
his literary talent carries over to the realm of logic, in particular
to the design of his logical exercises for students and his various
puzzles. So, it appears acceptable to use his literary pseudonym `Lewis
Carroll' throughout to refer to him.

Starting point is our encounter with the complete literary works edition
\cite{LC-Complete}, which, due to their literary qualities, nevertheless
contains also a number of exercises for his students of classical
logic and some other riddles. The opening exercise \cite[p. 1119, problem 1]{LC-Complete}
did serve as an early trigger to pursue a more thorough analysis of
his approach to classical logic, which we only now, decades later,
found the time to undertake. The exercise in question goes like this:

Given the assumptions

(1) Babies are illogical,

(2) Nobody is despised who can manage\emph{ }a crocodile, 

(3) Illogical persons are despised, 

what is the conclusion from this? The answer is given as: `Babies
cannot manage crocodiles'. As amusing as this conclusion sounds,
there are immediate questions arising from this, such as `Are babies
really illogical?', `Are illogical persons always despised?', `Is
really nobody despised just for being able to manage crocodiles?'.
The point is, that from a logical perspective, these concerns do not
matter. We are to assume that (1)-(4) is true and we are merely being
asked to find out what follows from this assumption. For illustrating
this point let us consider the completely nonsensical assumptions

(1) All cats understand French, 

(2) Some chicken are cats.

Then, it can be concluded that `Some chicken understand French'.
It's a standard so-called syllogism of a form already analysed by
Aristoteles. If we would have studied Aristotelian logic, then we
would know what \emph{`the conclusion'} means. It is apparently
\emph{not} just everything that follows. In a more general setting
it is not so obvious what \emph{`the conclusion'} could mean, and
part of the following is dedicated to giving a reasonable meaning
to this term.

On the way, we have to deal with a somewhat tricky situation. To understand
logic we already need to know a basic form of logic, otherwise we
will not be able to understand the logical arguments. Also, we need
to clarify, what we want to talk about and so we have to deal with
`things' and collections thereof and the notion of belonging to
this or that collection. These in turn are concepts properly dealt
with in the context of `set theory'. So, we need to have already
an understanding of sets to understand the objects of logic. Conversely
deepening our understanding of sets requires of course a better understanding
of logic. In the light of Lewis Carroll's approach to consider logic,
we encounter a third issue: algebra. His perspective is to think of
logic as a system of equations involving variables, which are only
allowed to be $1$, if the variable is supposed to represent a logical
statement that we consider to be `true', or $0$, if the variable
is supposed to represent a logical statement that we consider to be
`false'. If the variables stand for functions assigning truth-values,
i.e. $0$ or $1$, to something, we shall speak of predicates. This
point of view amounts to considering the arithmetic rules $a\cdot b=1$
only if $a=b=1$, otherwise $a\cdot b=0$, and $a+b=1$, if either
$a=1$ and $b=0$ or $a=0$ and $b=1$. In all other cases, i.e. if
$a=b$ we have $a+b=0$. Also, sets are better understood in terms
of their so-called characteristic functions, which give the value
$1$ if applied to something that belongs to a set, and $0$ if not.
So, this connects us to the arithmetic of $0$ and $1$. If we think
of $0$ standing for `even' and $1$ for `odd', then these rules
should be well-known from school days. It is just the arithmetic of
even and odd natural numbers, see subsection \ref{subsec:A-brief-look}.
Conversely, logic and sets are needed to better understand all these
statements about algebraic concepts. 

We call this situation the \emph{entangled trinity,} consisting of
algebraic operations in a setting based on arithmetic using only $0$
and $1$ (Boolean algebra), the concept of sets (set theory), and
logical arguments (predicate logic). Fortunately, in these enlightened
days, we may assume that the interested reader is somewhat familiar
with these basic concepts as a starting point of our journey into
a slightly deeper understanding of these topics.

The German mathematical logician Paul Lorenzen likened our situation
with somebody floating in the toiling sea while collecting driftwood
to improve his ability to keep afloat, eventually maybe being able
to build a raft. Our approach is to lift these three areas all at
the same time to a slightly higher level.

We shall speak of \emph{three facets} of the entangled trinity, each
eventually offering an equivalent representation. Indeed, somewhat
luckily, the discussion of $\left\{ 0,1\right\} $-valued (or Boolean-valued)
functions, sets and predicates eventually turn out to be the same
thing, which is the main benefit of this observation. 

Approaching logical issues with arithmetic means is really at the
core of Lewis Carroll's approach to classical logic, which in turn
is close to Boole's perspective, \cite{Boole1854}. His criticism
of Boolean logic is more from a didactical and intuitive perspective.
Our discussion actually consoles Dodgson's with Boole's approach showing
that in essence they are indeed the same, although looking quite different
in interpretation and usage. As for the methods Lewis Carroll suggested to solve his puzzles, we have to say that they do give a conclusion, but only in simple cases they yield the sharpest conclusion, which is indeed what we shall be focusing on, since weaker conclusions can then derived trivially. Also, we have restricted our references merely to those actually used in the work presented here.  

An obvious and actually intentional side-benefit of the Boolean arithmetic
approach is that computers are based on arithmetic with $0$ and $1$,
and so, so with $0$ for `false' and $1$ for `true', can conveniently
be used to explore logical questions. Let us start with a brief reminder
of the more familiar arithmetic with natural numbers.

\subsection{A brief look at the arithmetic of even and odd natural numbers}\label{subsec:A-brief-look}

Natural numbers are at the root of the constructions of many classical
mathematical concepts\footnote{For starters, we note that the power set $2^{\mathbb{N}}$ of the
natural numbers $\mathbb{N}$, i.e. the set of all subsets of $\mathbb{N}$,
has the same cardinality, that is the same (infinite) ``number of
elements'', as the set $\mathbb{R}$ of real numbers.} and that is also where we shall start with our journey. There is
a rather clever construction of natural numbers given by John von
Neumann (in the context of constructing higher ordinals) within the
language of elementary set theory, which goes recursively like this:
\begin{enumerate}
\item The empty set $\left\{ \right\} $, the set `with nothing in it',
is taken as representing $0$,
\item all larger numbers are given by the recursion 
\[
n+1\coloneqq n\cup\left\{ n\right\} ,
\]
i.e. by collecting everything in the already constructed set $n$
and adjoin it by this set itself.
\item the set of all so constructed sets is denoted by $\mathbb{N}$, the
set of natural numbers, the smallest infinite set\footnote{We note that, by definition, $0$ is a natural number
\[
0\in\mathbb{N},
\]
indeed the smallest one there is.}.
\end{enumerate}
From this natural numbers construction, the usual arithmetic of natural
numbers, i.e. the elements of $\mathbb{N}$, can be derived. As a
byproduct of von Neumann's construction, we have that
\[
n^{\prime}\coloneqq n+1=\left\{ 0,\ldots,n\right\} \text{ for natural numbers}\:n.
\]
The $n^{\prime}$ is called the successor of $n$ and is an often
used alternative name for $n+1$ to distinguish the initial construction
of $\mathbb{N}$ from the operation of addition, which is based on
the successor construction. From the knowledge of what $n^{\prime}$
is, we construct general sums of two numbers via
\begin{align*}
n+m^{\prime} & \coloneqq\left(n+m\right)^{\prime}.
\end{align*}
Similarly from
\[
n\cdot0=0\cdot n=0,
\]
we get recursively 
\[
n\cdot m^{\prime}=n\cdot m+n.
\]
Limited to the few tools available there is little more than `complete
induction' to then prove things like commutativity of multiplication,
i.e.
\[
n\cdot m-m\cdot n=0.
\]
For `$n$ is a natural number' we shall frequently write $n\in\mathbb{N}$,
(read as `$n$ element of $\mathbb{N}$'). As stated, von Neumann's
construction is based on set language, which displays already two
components of what we suggested to call the `entangled trinity',
namely set theory and predicate logic. 

As a curious\footnote{If we introduce strict order $<$ by writing $n<m$ if and only if
$n\leq m$ and $n\neq m$, then we find that $n<m$ is in von Neumann's
construction indeed the same as stating that $n$ is an element of
the set $m$, i.e. as a formula
\[
n\in m\,.
\]
} consequence of the von Neumann construction we find in particular
that for natural numbers $n,m,$ 
\[
n\leq m
\]
holds if and only if $n$ is a subset of $m$, i.e. as a formula
\[
n\subseteq m.
\]
In the form of a Venn diagram:
\begin{center}
\includegraphics{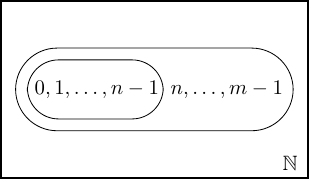}
\par\end{center}

The arithmetic we need for our purposes is fortunately much more simple
then the standard one in $\mathbb{N}$. Indeed, it is induced by the
standard arithmetic in $\mathbb{N}$ if we focus on the property of
natural numbers being divisible by $2$, even numbers, or not, odd
numbers. The smallest even number is $0$ and the smallest odd number
is $1$. Using $0$ and $1$ as labels for these two types of numbers,
we arrive at the rules
\begin{align*}
0+0 & =0,\\
0^{\prime}=0+1 & =1+0=1,\\
1^{\prime}=1+1 & =0,
\end{align*}
since, as we recall, $"even"+"even"="even"="odd"+"odd"$ and $"even"+"odd"="odd"="odd"+"even"$.
We also note that the successor operation has become circular
\begin{align*}
0^{\prime} & =1,\\
1^{\prime} & =0.
\end{align*}
After all, we only have even $"even"$ or $"odd"$ available. For
the multiplication we have
\begin{align*}
0\,\cdot\,0 & =0,\\
0\,\cdot\,1 & =1\,\cdot\,0=0,\\
1\,\cdot\,1 & =1,
\end{align*}
since $"even"\,\cdot\,"even"="even"="odd"\,\cdot\,"even"="even"\,\cdot\,"odd"$
and $"odd"\,\cdot\,"odd"="odd".$ This is the entry point of the first
facet of the `entangled trinity': arithmetic in the set $\left\{ 0,1\right\} .$

Picking up the von Neumann idea of constructing natural numbers, we
have in particular $2=\left\{ 0,1\right\} $, so that, using $2$
to abbreviate $\left\{ 0,1\right\} $, we have established the operations
``$+$'' and ``$\:\cdot\:$'' in the set $2$, indeed making $2$,
equipped with the operations ``$+$'' and ``$\:\cdot\:$'', the
smallest field possible. In particular, we have that $1$ is its own
additive and multiplicative inverse, since
\[
1+1=0,1\cdot1=1.
\]

We also have an obvious order inherited from the natural numbers:
\[
0\leq0,0\leq1,1\leq1.
\]
Note, however, that, in contrast to the situation for the order of
the natural numbers, we have that if 
\[
a\leq b
\]
 then\footnote{This shows that $2$ is not only a field, but also the smallest non-trivial
lattice, but not an ordered field in the strict sense, since the order
we have is not compatible with addition. In the natural numbers we
have this compatibility, since $a\leq b$ yields $a^{\prime}\leq b^{\prime}$.
Here in $2$ we have the reverse conclusion!} 
\[
b^{\prime}=1+b\leq1+a=a^{\prime}.
\]
We also note that $\leq$ can be expressed again in 2-arithmetic terms.
Indeed, we find that $a\leq b$ is the same as the equation
\begin{equation}
a\cdot b^{\prime}=0\label{eq:leq}
\end{equation}
or 
\[
a=a\cdot b.
\]
In form of a table we have
\begin{center}
\begin{tabular}{cc|ccc}
$a$ & $b$ & $a\leq b$ & $a\cdot b'=0$ & $a=a\cdot b$\tabularnewline
\hline 
0 & 0 & $0\leq0$ & $0\cdot0'=0$ & $0=0\cdot0$\tabularnewline
0 & 1 & $0\leq1$ & $0\cdot1'=0$ & $0=0\cdot1$\tabularnewline
1 & 1 & $1\leq1$ & $1\cdot1'=0$ & $1=1\cdot1$\tabularnewline
\hline 
1 & 0 & $1\not\leq0$ & $1\cdot0'\neq0$ & $1\neq1\cdot0$\tabularnewline
\end{tabular}
\par\end{center}

As a 2-arithmetic notational simplification, we shall henceforth
mostly drop the multiplication sign. This is normal practise in arithmetic
and Lewis Carroll follows the same convention. Thus, \eqref{eq:leq}
becomes
\begin{equation}
ab^{\prime}=0\label{eq:leq-1}
\end{equation}
or 
\[
a=ab.
\]

The rather simple setting developed so far is indeed, as we shall
see, the realm of logic -- entangled with sets and 2-arithmetic --
and so of the structure \emph{behind} the considerations in Lewis
Carroll's work on logic. This perspective explains a part of the title,
since rather than going \emph{through} the looking glass we shall
aim to peek \emph{behind} it.

With these preliminary remarks we should be ready to investigate our
topic somewhat more closely.

\section{On Function Algebra over $2$}

\subsection{Arithmetic in $2$ (continued)}

As we noted, the set 2 is the smallest (finite) field with a natural
order. It turns out, that instead of using $"+"$ and $"\:\cdot\:"$,
we really only need a single operation, the Sheffer stroke $"\:|\:"$
(in computer science and engineering this is relabelled as `NAND
gate'), given by
\begin{align*}
0|0 & =1,\\
0|1 & =1|0=1,\\
1|1 & =0,
\end{align*}
which is connected with the arithmetic operations via
\[
a|b=1+ab,\:a,b\in2.
\]
Conversely, we can express the arithmetic operations in terms of the
Sheffer stroke as
\begin{align*}
ab & =\left(a|b\right)|\left(a|b\right),\\
a+b & =\left(a|\left(b|b\right)\right)|\left(b|\left(a|a\right)\right),\:a,b\in2.
\end{align*}
In particular, we have 
\begin{align*}
a^{\prime} & =\left(1+a\right)=\left(a|a\right),\:a\in2.
\end{align*}
Nevertheless, it is frequently easier and shorter to understand an
expression using the more intuitive operators $"+"$ and $"\:\cdot\:"$
operations. Indeed, there are other operations, which facilitate the
use of our intuition, such as $"\max$'' and $"\min$''. They can
be defined in various notations as
\[
a\wedge b\coloneqq\wedge\left\{ a,b\right\} \coloneqq\min\left\{ a,b\right\} \coloneqq ab,
\]
\begin{align*}
a\vee b\coloneqq\vee\left\{ a,b\right\} \coloneqq\max\left\{ a,b\right\}  & \coloneqq a+b+ab,\:a,b\in2,
\end{align*}
and provide the formal presentation of what we would expect the maximum
and minimum of the two numbers $0$ and $1$ to be. The operation
symbols $\wedge,\vee$ are referred to as `meet' and `join'.

We note that
\[
\max\left\{ a,b\right\} +\min\left\{ a,b\right\} =a+b,\:a,b\in2.
\]
 In $2$ we get by a simple arithmetic calculation that
\begin{equation}
1+\max\left\{ a,b\right\} =\min\left\{ 1+a,1+b\right\} ,\:a,b\in2.\label{eq:maxmin}
\end{equation}

There is another single operation, the so-called Peirce arrow $"\:\downarrow\:"$
(renamed as `NOR gate' in engineering and computer science), which
also allows to express all other operations. It is given by
\begin{align*}
0\downarrow0 & =1,\\
0\downarrow1 & =1\downarrow0=0,\\
1\downarrow1 & =0,
\end{align*}
which is connected with the arithmetic operations via
\[
a\downarrow b=1+\max\left\{ a,b\right\} ,\:a,b\in2.
\]
Conversely, we can also express the arithmetic operations in terms
of the Peirce arrow as
\begin{align*}
a\cdot b & =\left(a\downarrow a\right)\downarrow\left(b\downarrow b\right),\\
a+b & =\left(a\downarrow b\right)\downarrow\left(\left(a\downarrow a\right)\downarrow\left(b\downarrow b\right)\right),\:a,b\in2.
\end{align*}
We also note that
\[
a^{\prime}=a\downarrow a,\:a\in2.
\]
For convenience we may also introduce $\max$ and $\min$ for larger
(finite) sets recursively 
\begin{align}
\bigvee\left\{ a_{0}\right\} \coloneqq\max\left\{ a_{0}\right\}  & =\bigwedge\left\{ a_{0}\right\} \coloneqq\min\left\{ a_{0}\right\} \coloneqq a_{0}\nonumber \\
\bigvee\left\{ a_{0},a_{1},\ldots,a_{n}\right\} &\coloneqq\max\left\{ a_{0},a_{1},\ldots,a_{n}\right\}\notag\\  & \coloneqq\max\left\{ \max\left\{ a_{0},a_{1},\ldots,a_{n-1}\right\} ,a_{n}\right\} \label{eq:bigvee-bigwedge}\\
\bigwedge\left\{ a_{0},a_{1},\ldots,a_{n}\right\} &\coloneqq\min\left\{ a_{0},a_{1},\ldots,a_{n}\right\}\notag\\  & \coloneqq\min\left\{ \min\left\{ a_{0},a_{1},\ldots,a_{n-1}\right\} ,a_{n}\right\} ,\nonumber \\
 & a_{k}\in2,\:k\in\mathbb{N}.\nonumber 
\end{align}
These notations for the `big' max and `big' min go back to Garrett
Birkhoff introduced in connection with lattice theory, from where
they diffused into logic.

With this we can generalise \eqref{eq:maxmin} for $a_{k}\in2,\:k,n\in\mathbb{N}$ to
\begin{equation}
1+\max\left\{ a_{0},a_{1},\ldots,a_{n}\right\} =\min\left\{ 1+a_{0},1+a_{1},\ldots,1+a_{n}\right\}.\label{eq:MaxMin}
\end{equation}
A literal is a letter standing for a value in $2$ or its complement.
Thus, if $a$ is a literal so is $1+a$. In case of $\max$-of-literals
expression or a $\min$-of-literals expression one occasionally speaks
concisely of maxterm and minterm, respectively. 

We finally note that even $"="$ can be expressed in 2-arithmetic
terms. Indeed,
\[
\left(a=b\right)=\left(1+a+b\right),
\]
which is, however, obviously circular, but still of some interest.

Since
\[
\min\left\{ a,\max\left\{ b,c\right\} \right\} =\max\left\{ \min\left\{ a,b\right\} ,\min\left\{ a,c\right\} \right\} ,
\]
\[
\max\left\{ a,\min\left\{ b,c\right\} \right\} =\min\left\{ \max\left\{ a,b\right\} ,\max\left\{ a,c\right\} \right\} ,
\]
\begin{align*}
\left(\max\left\{ a,b\right\} \right)^{\prime} & =\min\left\{ a^{\prime},b^{\prime}\right\} ,\\
\left(\min\left\{ a,b\right\} \right)^{\prime} & =\max\left\{ a^{\prime},b^{\prime}\right\} ,
\end{align*}
every arithmetic expressions can be written as a maxterm of minterms
(a normal form with negation). For example, we have

\begin{align*}
\left(ab\right) & =\min\left\{ a,b\right\} =\max\left\{ \min\left\{ a,b\right\} \right\} ,
\end{align*}
\begin{align*}
\left(a+b\right) & =\left(\max\left\{ a,b\right\} \left(\min\left\{ a,b\right\} \right)^{\prime}\right),\\
 & =\left(\max\left\{ a,b\right\} \,\cdot\,\max\left\{ a^{\prime},b^{\prime}\right\} \right),\\
 & =\left(\min\left\{ \max\left\{ a,b\right\} ,\max\left\{ a^{\prime},b^{\prime}\right\} \right\} \right),\\
 & =\left(\max\left\{ \min\left\{ a,b^{\prime}\right\} ,\min\left\{ b,a^{\prime}\right\} \right\} \right),
\end{align*}
\begin{align*}
\left(a\leq b\right) & =\left(\min\left\{ a,b^{\prime}\right\} \right)^{\prime}=\max\left\{ a^{\prime},b\right\} ,
\end{align*}
\begin{align*}
\left(a=b\right) & =\left(\max\left\{ \min\left\{ a,b^{\prime}\right\} ,\min\left\{ b,a^{\prime}\right\} \right\} \right)^{\prime},\\
 & =\min\left\{ \max\left\{ a,b^{\prime}\right\} ,\max\left\{ b,a^{\prime}\right\} \right\} ,\\
 & =\max\left\{ \min\left\{ a,b\right\} ,\min\left\{ a^{\prime},b^{\prime}\right\} \right\} .
\end{align*}

Written as nullities in maxterm of minterms normal form, we have for
this
\begin{align*}
\left(ab\right) & =\left(\left(\min\left\{ a,b\right\} \right)^{\prime}=0\right)\\
 & =\left(\max\left\{ a^{\prime},b^{\prime}\right\} =0\right)
\end{align*}
\begin{align*}
\left(a+b\right) & =\left(\left(\min\left\{ \max\left\{ a,b\right\} ,\max\left\{ a^{\prime},b^{\prime}\right\} \right\} \right)^{\prime}=0\right)\\
 & =\left(\max\left\{ \min\left\{ a,b\right\} ,\min\left\{ a^{\prime},b^{\prime}\right\} \right\} =0\right)
\end{align*}
\begin{align*}
\left(a\leq b\right) & =\left(\max\left\{ \min\left\{ a,b^{\prime}\right\} \right\} =\min\left\{ a,b^{\prime}\right\} =0\right)
\end{align*}
\begin{align*}
\left(a\neq b\right) & =\max\left\{ \min\left\{ a,b^{\prime}\right\} ,\min\left\{ b,a^{\prime}\right\} \right\} \\
 & =\left(\left(\max\left\{ \min\left\{ a,b^{\prime}\right\} ,\min\left\{ b,a^{\prime}\right\} \right\} \right)^{\prime}=0\right)\\
 & =\left(\min\left\{ \max\left\{ a,b^{\prime}\right\} ,\max\left\{ b,a^{\prime}\right\} \right\} =0\right)\\
 & =\left(\max\left\{ \min\left\{ a,b\right\} ,\min\left\{ b^{\prime},a^{\prime}\right\} \right\} =0\right)
\end{align*}

In general, any $2$-arithmetic term involving a set of names $\left\{ a_{0},\ldots,a_{N-1}\right\} $
standing for possible values in $2$ can be written as 
\[
\max\left\{ \min A_{0},\ldots,\min A_{M-1}\right\} ,
\]
where $A_{k},\:k=0,\ldots,M-1$, stands for certain subsets of $\{ a_{0},\ldots,a_{N-1},$ $a_{0}^{\prime},\ldots,a_{N-1}^{\prime}\}$.
This will be the standard reference form, a normal form, which we
will base our analysis on.

In terms of the `meet' and `join' operations introduced above,
this becomes particularly clear, since they show that behind the interaction
of $\max$ and $\min$ utilised here there is a distributive law between
`meet' $\wedge$ and `join' $\vee$: 
\[
\left(a\wedge\left(b\vee c\right)\right)=\left(\left(a\wedge b\right)\vee\left(a\wedge c\right)\right),
\]
\[
\left(a\vee\left(b\wedge c\right)\right)=\left(\left(a\vee b\right)\wedge\left(a\vee c\right)\right).
\]
The interaction of $\max$ and $\min$ with the successor operation
appears now in terms of these operations as a consequence of the so-called
de Morgan rules
\begin{align*}
\left(a\vee b\right)^{\prime} & =a^{\prime}\wedge b^{\prime},\\
\left(a\wedge b\right)^{\prime} & =a^{\prime}\vee b^{\prime}.
\end{align*}

\subsection{Lewis Carroll's approach in arithmetic terms}\label{subsec:Lewis-Carroll's-approach}

Although, it may at this point look a bit disconnected from logic,
since the logic facet of the entangled trinity has not been built
up yet, this is already the logical place to initiate the discussion
of Lewis Carroll's approach, since his approach is indeed arithmetic
in nature.

Our starting point is an expression of the form
\[
\mathcal{A}=\max\left\{ \min A_{0},\ldots,\min A_{M-1}\right\} .
\]

On inspection of a particular sets $A_{k},\:k=0,\ldots,M-1$, we may
have some names in $\left\{ a_{0},\ldots,a_{N-1}\right\} $ with also
its successor appearing. These will be called eliminands, since we
will want to eliminate those from consideration. The remaining elements
are referred to as retinends, since these we shall want to keep. Note
that this classification depends on the particular form $\mathcal{A}$!
There are forms equivalent to $\mathcal{A}$, where eliminands and
retinends are different.

To make this more palpable let us look at a simple example. Consider
\begin{equation}
\mathcal{A}=\max\left\{ \min\left\{ a,b\right\} ,\min\left\{ b^{\prime},c\right\} \right\} =\max\left\{ ab,b^{\prime}c\right\} \label{eq:Ass}
\end{equation}
$b$ is the only eliminand, $a,c$ are retinends. Inspecting
\[
\tilde{\mathcal{A}}=\max\left\{ abc,abc^{\prime},b^{\prime}c\right\} 
\]
we see that the values of $\mathcal{A}$ and $\mathcal{\tilde{\ensuremath{A}}}$
for the same values of $a,b,c$ are the same\footnote{The sensitivity to the particular form of the underlying expression
is essential for Lewis Carroll's approach to logical problems and
it strictly speaking more of a linguistic topic than an inner logical
issue. This shows his closeness to classical logic rather than the
emerging mathematical logic of his time. His criticism of Boole's
approach to logic may well rest on this particular issue. Nevertheless,
after constructing the appropriate wording of the problem in logical
terms, the solution process remains purely mathematical in nature
and coalesces with Boole's perspective. }. But in $\mathcal{\tilde{\ensuremath{A}}}$ now $b$ and $c$ are
eliminands, $a$ is the only retinend. 

Our interpretation of Lewis Carroll's approach now is that he wants
to find the largest\footnote{We choose to focus on the largest lower bound since any other lower bound can easily be derived from the largest by further constraints. This approach is actually a sharpening of Lewis Carroll's approach.} lower bound in terms of retinends only! In the
first instance this is
\begin{equation}
\mathcal{A}\geq ac\eqqcolon\mathcal{C},\label{eq:Conc}
\end{equation}
in the second case the result is the uninteresting estimate
\[
\tilde{\mathcal{A}}\geq0.
\]
Lewis Carroll's perspective on
\[
\mathcal{A}\geq\mathcal{C}
\]
is to consider the equivalent form
\[
\mathcal{A}^{\prime}\leq\mathcal{C}^{\prime}.
\]
Noting that
\[
\mathcal{A}^{\prime}=\left(\mathcal{A}=0\right),\mathcal{C}^{\prime}=\left(\mathcal{C}=0\right),
\]
this says, if 
\[
\mathcal{A}=0
\]
 then 
\[
\mathcal{C}=0.
\]
Now $\mathcal{A}=\max\left\{ \min A_{0},\ldots,\min A_{M-1}\right\} =0$
is equivalent to the system of equations
\[
\min A_{0}=0,\ldots,\min A_{M-1}=0,
\]
which is, as we shall see later, is precisely Lewis Carroll's take
on how to consider logical puzzles. There is, however, still quite
a bit to go until we are able to address this more precisely.

\subsection{$2$-valued functions }

At least, we are now ready to complete the introduction of the first
facet of the entangled trinity. Since we are focusing on classical
aspect we will assume here that our `universe of discourse' $\mathcal{U}$
is finite throughout and of course $\mathcal{U}$ non-empty, otherwise
little -- or rather too much -- can be said about its content. Our
interest is now in 2-valued functions, which together with image-wise
arithmetic operations, we shall introduce explicitly here, yields
a structure known as an algebra (Boolean algebra). We shall refer
to this algebra briefly as a $2$-algebra. The functions $P$ in this
2-algebra assign to any given $x\in\mathcal{U}$ uniquely either $0$
or $1$, so that the result of $P$ being applied to $x\in\mathcal{U}$
results in a specific value in $2$ denoted by $P\left(x\right)$.
The collection of such functions will be denoted by $2^{\mathcal{U}}$,
which is the set of all functions with arguments ranging in $\mathcal{U}$
and having values in $2=\left\{ 0,1\right\} $. 

Let us explicitly lift the arithmetic in $2$ to a 2-algebra of functions
$P,Q,$ by defining image-wise operations:
\begin{align}
P^{\prime}\left(x\right) & \coloneqq\left(P\left(x\right)\right)^{\prime}\nonumber \\
\left(PQ\right)\left(x\right) & \coloneqq P\left(x\right)Q\left(x\right)\nonumber \\
\left(\min\left\{ P,Q\right\} \right)\left(x\right) & \coloneqq\min\left\{ P\left(x\right),Q\left(x\right)\right\} =P\left(x\right)\wedge Q\left(x\right)\notag\\ &=P\left(x\right)Q\left(x\right)\eqqcolon\left(P\wedge Q\right)\left(x\right)\nonumber \\
\left(P+Q\right)\left(x\right) & \coloneqq P\left(x\right)+Q\left(x\right)\label{eq:lift}\\
\left(\max\left\{ P,Q\right\} \right)\left(x\right) & \coloneqq\max\left\{ P\left(x\right),Q\left(x\right)\right\} =P\left(x\right)\vee Q\left(x\right)\eqqcolon\left(P\vee Q\right)\left(x\right)\nonumber \\
\left(P\leq Q\right)\left(x\right) & \coloneqq\left(P\left(x\right)\leq Q\left(x\right)\right)\nonumber \\
\left(P=Q\right)\left(x\right) & \coloneqq\left(P\left(x\right)=Q\left(x\right)\right)\nonumber 
\end{align}
for all $x\in\mathcal{\ensuremath{U}}.$ Also the notation \eqref{eq:bigvee-bigwedge}
carries over to the 2-function setting via
\begin{align}
\left(\bigvee\left\{ P_{0},P_{1},\ldots,P_{n}\right\} \right)\left(x\right) & \coloneqq\bigvee\left\{ P_{0}\left(x\right),P_{1}\left(x\right),\ldots,P_{n}\left(x\right)\right\} ,\label{eq:bigvee-bigwedge-functions}\\
\left(\bigwedge\left\{ P_{0},P_{1},\ldots,P_{n}\right\} \right)\left(x\right) & \coloneqq\bigwedge\left\{ P_{0}\left(x\right),P_{1}\left(x\right),\ldots,P_{n}\left(x\right)\right\} ,\:x\in\mathcal{U}.\nonumber 
\end{align}
Frequently, these are written in an index notation analogous to the
big summation sign as
\[
\bigvee_{k\in n+1}P_{k}\coloneqq\bigvee_{k=0}^{n}P_{k}\coloneqq\bigvee\left\{ P_{0},P_{1},\ldots,P_{n}\right\} ,
\]
\[
\bigwedge_{k\in n+1}P_{k}\coloneqq\bigwedge_{k=0}^{n}P_{k}\coloneqq\bigwedge\left\{ P_{0},P_{1},\ldots,P_{n}\right\} ,\:n\in\mathbb{N}.
\]

\section{Sets}

We clearly have used a vague concepts of sets already, but now we
are prepared to dig a bit deeper in order to establish the second
facet of the entangled trinity. We shall identify sets with their
so called characteristic function, which gives 1 if the argument belongs
to the set and $0$ otherwise, by reusing the same name. This results
in a 2-valued function defining a set and the function is in turn
the \emph{characteristic function} of the set:
\[
\left\{ x|P\left(x\right)\right\} \coloneqq\left\{ x|P\left(x\right)=1\right\} \coloneqq P.
\]
The ``element of'' relation $\in$ is then re-linked via interpreting
$x\in P$ simply as meaning that $P\left(x\right)=1$, where $P$
is a 2-valued function on a given finite universe $\mathcal{U}$.
Thus, sets are essentially just elements of a 2-algebra over a fixed
domain $\mathcal{U}$. In particular, the earlier used (naive) set
of functions $2^{\mathcal{U}}$ is the same as the set of all subsets.
In the context of sets we also consider larger universes such as $2^{\mathcal{U}}$
or even higher power sets (sets of sets of ...).

We see the entanglement by noting that mappings, a concept we had
to use already, are of course set theoretical objects, which we are
only now addressing in more depth. 

We have earlier already used set notations of the form
\[
\left\{ a_{0},a_{1},\ldots,a_{n-1}\right\} ,
\]
which we now see as properly defined by noting that the empty set
is given by
\begin{align*}
\left\{ \right\}  & \coloneqq\left\{ x|\left(x\neq x\right)=1\right\} ,
\end{align*}
and a set with only one element, a singleton, is given as 
\begin{align*}
\left\{ a\right\}  & \coloneqq\left\{ x|\left(x=a\right)=1\right\} .
\end{align*}
Here we treat $x\neq x$ and $x=a$ as mappings into $2$. This is
an anticipation of the logic facet of the entangled trinity, since
we are using $0,1$ as `false' and `true' here already. So, since
$x\neq x$ is always `false' the value of $x\neq x$ can never be
$1$. similarly $x=a$ can only be true if $x$ is actually equal
to $a$. This appears to be quite circular, but it is the typical
situation of entanglement we have to deal with. We always must assume
that we know things on a more naive level already to be able to make
them more precise later. 

Then, recursively
\begin{align*}
\left\{ a_{0},a_{1},\ldots,a_{n}\right\}  & \coloneqq\left\{ a_{0},a_{1},\ldots,a_{n-1}\right\} \cup\left\{ a_{n}\right\} ,\,n\in\mathbb{N}.
\end{align*}
Likewise we can now in the framework of set theory make precise what
a mapping is, which is another term we have used already earlier on
in an intuitive sense.

Starting point for explaining what a mapping actually is, lies in
the concept of ordered pairs.

An ordered pair $\left(a,b\right)$ is by definition a set of subsets
of the form $\left\{ \left\{ a\right\}\!,\!\left\{ a,b\right\} \right\} $,
($a,b$ taken -- say -- from a specified universe $\mathcal{U}$,
which is the classical perspective being analysed here). If an ordered
pair $\left(a,b\right)\coloneqq\left\{ \left\{ a\right\} ,\left\{ a,b\right\} \right\} $
is given, then $a$ can be easily recovered. Indeed, 
\[
\left\{ a\right\} =\bigcap\left\{ \left\{ a\right\} ,\left\{ a,b\right\} \right\} 
\]
and $a$ is the only element\footnote{If we may assume that $a$ is itself a subset, that is -- by our
identification -- a characteristic function of a given universe,
then 
\[
a=\bigcap\left\{ a\right\} =\bigvee\left\{ a\right\} .
\]
} of $\left\{ a\right\} $. Now, if $b=a$ we have also recovered $b$.
If, however, $b\neq a$ then\footnote{Since $P_{a}\left(x\right)\coloneqq\left(x=a\right)$ defines the
characteristic function of $\left\{ a\right\} $, then $\left\{ a\right\} ^{\prime}$
is the set $\left\{ x|P_{a}^{\prime}\left(x\right)=1\right\} $, i.e.
$P_{a}^{\prime}$ is the characteristic function of $\left\{ a\right\} ^{\prime}$.} 
\begin{align*}
\left\{ b\right\}  & =\left\{ a\right\} ^{\prime}\cap\bigcup\left\{ \left\{ a\right\} ,\left\{ a,b\right\} \right\} \\
 & =\left\{ a,b\right\} \cap\left\{ a\right\} ^{\prime}
\end{align*}
and $b$ is recovered as the only element of $\left\{ b\right\} $.
As already stated, an ordered pair of elements in $\mathcal{U}$ is
actually a set of subsets of $\mathcal{U}$. Note that since the set
of subsets is $2^{\mathcal{U}}$, sets of such subsets are actually
in $2^{2^{\mathcal{U}}}=2^{\left(2^{\mathcal{U}}\right)}$. 

Also note that the order in which $a,b$ appear is crucial, hence
\emph{ordered} pair. Indeed, if we have 
\[
\left(a,b\right)=\left(b,a\right),
\]
which means 
\[
\left\{ \left\{ a\right\} ,\left\{ a,b\right\} \right\} =\left\{ \left\{ b\right\} ,\left\{ a,b\right\} \right\} ,
\]
we obtain
\[
\left\{ a\right\} =\bigcap\left\{ \left\{ a\right\} ,\left\{ a,b\right\} \right\} =\bigcap\left\{ \left\{ b\right\} ,\left\{ a,b\right\} \right\} =\left\{ b\right\} ,
\]
from which $a=b$ follows. 

A set containing only ordered pairs is called a relation. For example
the set
\[
\left\{ x|\left(x=\left(c,c\right)\right)\cdot\left(c\in\mathcal{U}\right)=1\right\} 
\]
describes the equality relation in $\mathcal{U}$. Here we have also
used $c\in\mathcal{U}$ as a mapping into $2$, where we assign\footnote{In other words, if we define $=$ as a 2-valued function on pairs
\[
\left(=\left(a,b\right)\right)\coloneqq\left(a=b\right),\:a,b\in\mathcal{U},
\]
then $=$ is the characteristic function of $\left\{ x|x=\left(c,c\right)\left(c\in\mathcal{U}\right)=1\right\} $.} $1$ as the value if (and only if) $c$ is an element of $\mathcal{U}$.
For emphasis, we have used here the multiplication sign, which we
normally agreed to leave out. Multiplication functions here as an
`and'. Once we get to introduce the logic facet of the entangled
trinity, this will become clearer. For sake of simplifying notation
we take a range of liberties to denote sets. So, instead of multiplying
the requirements for the elements of a set we may use just commas.
Also using the format $\left\{ x|P\left(x\right)=1\right\} $ may
some times be too cumbersome. We will also drop the `$=1$'-part\footnote{Once we have established the predicate logic facet of the entangled
trinity (see next section), we shall see that this convention is actually
quite reasonable.}. For example, we may write $\left\{ x|\left(x=\left(c,c\right)\right)\left(c\in\mathcal{U}\right)=1\right\} $
more compactly as
\[
\left\{ \left(c,c\right)|c\in\mathcal{U}\right\} .
\]

A particular set of ordered pairs is, if we take all ordered pairs
$\left(a,b\right)$ with $a\in M$ and $b\in N$, where $M,N$ are
arbitrary sets. This is known as the Cartesian product of $M$ and
$N$, given in our more liberal set notation as
\[
M\times N\coloneqq\left\{ \left(a,b\right)|a\in M,b\in N\right\} 
\]

A general relation $R$ is called a function if it is `right unique'.
This means if $a$ is given then there is at most one pair $\left(a,b\right)$
in $R$. 

One frequently wants to also record, where $a$ and $b,c$ are taken
from in the pair building process. In other words, from which Cartesian
product $M\times N$ the function $R$ is a subset of. This is noted
in a slightly overbearing form to invoke the intuitive image of mapping
something from some set into another set, as $R:\dom\left(R\right)\subseteq M\to N$
and is referred to as a mapping. Here 
\[
\dom\left(R\right)\coloneqq\left\{ a|\left(a,b\right)\in R\right\} 
\]
is called the domain of the mapping $R$. If the domain is all of
$M$ we leave out the `$\dom\left(R\right)\subseteq$'-part. The
analogue
\[
\ran\left(R\right)\coloneqq\left\{ b|\left(a,b\right)\in R\right\} 
\]
 is called the range of the mapping $R$. Sometimes it is preferred
to write even more detailed and indeed more symmetrically
\[
R:\dom\left(R\right)\subseteq M\to\ran\left(R\right)\subseteq N
\]
for the mapping $R$.

Since here $R$ is a function, we have for every $a\in\dom\left(R\right)$
a unique element $R\left(a\right)$ of $\ran\left(R\right)\subseteq N$
such that $\left(a,R\left(a\right)\right)\in R$, which may also be
detailed as the following `maps to' statement 
\[
a\mapsto R\left(a\right).
\]
Frequently $R$ is given by an explicit expression. For example, we
may consider $\subseteq$ as a mapping
\[
\subseteq:2^{\mathcal{U}}\times2^{\mathcal{U}}\to2,
\]
where
\[
\left(A,B\right)\mapsto\left(A\leq B\right).
\]
Recall that $2^{\mathcal{U}}$ denotes the so-called power set of
$\mathcal{U}$, which is the set of all subsets of $\mathcal{U}$,
and the set of all mappings from $\mathcal{U}$ to $2$. 

Generally, the set of mappings from $M$ to $N$ with $\dom\left(R\right)=M$
is denoted by
\[
N^{M},
\]
which is motivated by the observation that if $\#M\in\mathbb{N}$
and $\#N\in\mathbb{N}$ are the numbers of elements in $M$ and $N$,
respectively, then the number of mappings mapping from all of $M$
to $N$ is given by the number
\[
\#N^{\#M}.
\]
Despite $N^{M}$ being a `power' notation, the name power set is
reserved for sets of the form $2^{\mathcal{U}}$. We now understand
why $2^{\mathcal{U}}$ has been used to describe the set of all subsets.
By identifying sets with their characteristic functions, we have indeed
that the set of all subsets of $\mathcal{U}$ coincides with the set
of characteristic functions, which is $2^{\mathcal{U}}$.

We note here that the concept of a Cartesian product can be expanded
to any number of factors. Indeed, let $N_{k}\subseteq N,\:k=0,\ldots,M-1,$
$M\in\mathbb{N}$, then 
\[
N_{0}\times\cdots\times N_{M-1}\coloneqq\left\{ f:M\to N|f\left(k\right)\in N_{k},k=0,\ldots,M-1\right\} 
\]
denotes their Cartesian product, as a subset of $N^{M}$. The elements
of this Cartesian product are frequently noted by writing down a list
of the values parametrised by by $M$, i.e. for $f\in N_{0}\times\cdots\times N_{M-1}$
we may write
\[
\left(f\left(k\right)\right)_{k\in M}
\]
or 
\[
\left(f\left(k\right)\right)_{k=0,\ldots,M-1}
\]
or even
\[
\left(f\left(0\right),\ldots,f\left(M-1\right)\right).
\]
These notations are referred to as $M$-tuples and the entries $f\left(k\right)$, for $k=0,\ldots,M-1$,
are called components of the $M$-tupel. 

So, for $M=2$ we recover the ordered pair ($2$-tupel) and the initially
introduced Cartesian product $N_{0}\times N_{1}$. This shows that
$N_{0}\times\cdots\times N_{M-1}$ is indeed a generalisation of the
2-factor Cartesian product we started out with. If the factors are
all the same, i.e. $N=N_{k}$, $k=0,\ldots,M-1$, then we have
\[
N^{M}=N_{0}\times\cdots\times N_{M-1},
\]
which neatly supports the power notation for mappings from $M$ to
$N$.

As it is customary to use different notations in set theory than in
2-functions algebra, there is an emerging ``doublespeak'' motivated
by our identification of sets and $2$-valued functions. We have `union'
$\cup$, `intersection' $\cap$, and `subset' $\subseteq$ as
new terms for things we already know in the context of $2$-algebras.
We have 
\[
\left(A\cap B\right)\coloneqq\min\left\{ A,B\right\} =\left(A\wedge B\right)=\left(AB\right)
\]
\[
\left(A\cup B\right)\coloneqq\left(A\vee B\right)=\max\left\{ A,B\right\} 
\]
\begin{align*}
\left(A\subseteq B\right) & \coloneqq\left(A\leq B\right)\\
\left(A=B\right) & \coloneqq\left(A=B\right).
\end{align*}
The last equality is referring to the equality of functions, which
we already identified with the equality of sets at the outset, hence
the same notation. The `doublespeak' also is reflected in the natural
language, where we say something like `$A$ intersected with $B$'
or shorter `$A$ intersection $B$' instead of the earlier `$A$
times $B$'. For `max$\left\{ A,B\right\} $' we now may also speak
of `the union of $A$ and $B$' or `$A$ union $B$', for short.
Similarly, we say `$A$ is contained in $B$' or `all elements
of $A$ are in $B$' or briefly `all $A$ are $B$' instead of
the earlier `$A\leq B$ `. A notation that is unchanged is the successor
notation, but in the context of sets $A^{\prime}$ is not considered
as the lifted successor operation but as taking the set of everything
in the universe $\mathcal{U}$, which is not in $A$ and is referred
to as the `complement' of $A$.

The empty set $0=\left\{ \quad\right\} $ is identified with the constant
function $0$, the function always giving the value $0$, whereas
$\mathcal{U}$ is identified with the constant function $1$ always
giving the value $1$. 

Finally, we record a lesser used set notation, the ``symmetric difference''
\[
P\triangle Q\coloneqq\left(P\cap Q^{\prime}\right)\cup\left(P^{\prime}\cap Q\right),
\]
which is simply the same as $P+Q$ in terms of the function algebra
$2^{\mathcal{\ensuremath{U}}}$. The notation of \eqref{eq:bigvee-bigwedge-functions}
is now rewritten as a `big union' and a `big intersection':
\begin{align}
\left(\bigcup\left\{ P_{0},P_{1},\ldots,P_{n+1}\right\} \right) & \coloneqq\left(\bigvee\left\{ P_{0},P_{1},\ldots,P_{n+1}\right\} \right),\label{eq:bigcup-bigcap}\\
\left(\bigcap\left\{ P_{0},P_{1},\ldots,P_{n+1}\right\} \right) & \coloneqq\left(\bigwedge\left\{ P_{0},P_{1},\ldots,P_{n+1}\right\} \right).\nonumber 
\end{align}
These also are frequently written in an index notation as
\[
\bigcup_{k\in n+1}P_{k}\coloneqq\bigcup_{k=0}^{n}P_{k}\coloneqq\bigcup\left\{ P_{0},P_{1},\ldots,P_{n}\right\} ,
\]
\[
\bigcap_{k\in n+1}P_{k}\coloneqq\bigcap_{k=0}^{n}P_{k}\coloneqq\bigcap\left\{ P_{0},P_{1},\ldots,P_{n}\right\} ,\:n\in\mathbb{N}.
\]

Clearly, the function algebra $2^{\mathcal{\ensuremath{U}}}$ and
set theory are inseparably entangled in as much as they are all already
used when trying to explain what they are in more precise terms. Our
discussion of the function algebra $2^{\mathcal{U}}$ presupposed
at least a proto-concept of sets (in the form of a meta-language use),
and sets have been constructed here using the function algebra of
their characteristic functions. Indeed, neither works without the
other. Also, by our choice of starting point the other two facets
now appear as merely a repetition with a changed terminology. However,
we should keep in mind that we could have started with any of the
other facets, making then the rest look as doing the same with a different
terminology. Anyway, in our chosen approach we are now ready to look
at the third and final facet of the entangled trinity, which is predicate
logic.

\section{Predicate Logic (First Order)}

\subsection{Propositional logic }

Predicate logic is based on propositional logic in the same way as
$2$-algebra is based on $2$-arithmetic. We are now at the point
of introducing a `triplespeak' to describe the logic perspective
of our topic. Indeed, there is little more to do than just making
that transition to another language about essentially the same topic.
Indeed, looking at the basic 2-arithmetic, we now interpret 0 and
1, as `false' and `true', respectively, and think of $a,b,c,\cdots$
as statements of which `truth' and `falsehood' can be claimed,
i.e. propositions. This interpretation gives us access to structural
considerations about natural language statements, which is indeed
our final goal in all this. We need to be aware of the pitfalls coming
from the fact that there are syntactically and grammatically correct
statements, which are meaningless, so that the truth values `true',
$1$, or `false', $0$, cannot be assigned, and so they are not
propositions. In this context we call $"\wedge"$ `and' instead
of `meet' and $"\vee"$ `or' instead of `join'. We call $a^{\prime}$
the `negation' of $a$ or `not $a$'. There is commonly even a
different notation used for `not $a$', which is $\neg a$, but
of course $\neg a$ is the same as $a^{\prime}$. There is also implication
$\Rightarrow$, and equivalence $\Leftrightarrow$. We say $a$ implies
$b$ if $a\leq b$. Indeed, 
\begin{align*}
\left(0\Rightarrow0\right) & =\left(0\leq0\right)=1,\\
\left(0\Rightarrow1\right) & =\left(0\leq1\right)=1,\\
\left(1\Rightarrow0\right) & =\left(1\leq0\right)=0,\\
\left(1\Rightarrow1\right) & =\left(1\leq1\right)=1,
\end{align*}
and 
\begin{align*}
\left(0\Leftrightarrow0\right) & =\left(0=0\right)=1,\\
\left(0\Leftrightarrow1\right) & =\left(0=1\right)=0,\\
\left(1\Leftrightarrow0\right) & =\left(1=0\right)=0,\\
\left(1\Leftrightarrow1\right) & =\left(1=1\right)=1.
\end{align*}
The exclusive or\footnote{One also finds $\neq$ and -- unfortunately more common -- $\oplus$
for `exclusive or'. The latter already has an established meaning
in too many areas as a so-called `direct sum', which collides with
the use as `exclusive or'.} $\underline{\vee}$ is given by
\begin{align*}
\left(0\underline{\vee}0\right) & =0=\left(0+0\right),\\
\left(1\underline{\vee}0\right) & =1=\left(1+0\right),\\
\left(0\underline{\vee}1\right) & =1=\left(0+1\right),\\
\left(1\underline{\vee}1\right) & =0=\left(1+1\right).
\end{align*}
We see that by our identification we have that 
\[
+,\underline{\vee}
\]
are the same and also equal to $\neq$ as well as the negation $\nLeftrightarrow$
of $\Leftrightarrow$. 

The earlier mentioned Sheffer stroke and Peirce arrow are originally
inventions in the context of logic rather than 2-algebra, where we
chose to introduce them first.

\subsection{First Order Predicate logic}

This is the third and -- in our presentation -- the final facet
of the entangled trinity. First order predicate logic is just another
perspective of the 2-algebra (Boolean algebra) we introduced initially.
The 2-valued functions defined on $\mathcal{U}$ are now called predicates
in the context of logic and are interpreted as proposition about a
variable allowed to range in a fixed (finite!) universe of discourse
$\mathcal{U}$. Since in effect we identify characteristic functions
of subsets of a universe $\mathcal{U}$ with predicates defined on
$\mathcal{U}$, we have for a set $P=\left\{ x|P\left(x\right)=1\right\} $
that also
\[
P=\left\{ x|P\left(x\right)\Leftrightarrow1\right\} =\left\{ x|P\left(x\right)\right\} .
\]
Indeed, 
\[
P\left(x\right)=\left(P\left(x\right)\Leftrightarrow1\right),
\]
as one easily checks. In other words, $P$ is the set of all $x$
such that $P\left(x\right)$ is true. The two newly introduced notations
$\Rightarrow$ and $\Leftrightarrow$ are defined by
\begin{align*}
\left(P\Rightarrow Q\right) & \coloneqq\left(P\leq Q\right),\\
\left(P\Leftrightarrow Q\right) & \coloneqq\left(P=Q\right).
\end{align*}
These are both implicit all-statements with respect to the arguments
from the universe $\mathcal{U}$. 

In the context of logic $\bigvee$ and $\bigwedge$ reappear also
in a slightly different flavour as existence and all quantifiers.
In their simplest form they are defined\footnote{Note that the set $\left\{ P\left(x\right)|x\in\mathcal{U}\right\} $
can only be either $\left\{ 0\right\} $, or $\left\{ 1\right\} $,
or $\left\{ 0,1\right\} $, since sets are called equal if they have
the same elements.} by 
\begin{align*}
\bigvee_{x\in\mathcal{U}}P\left(x\right) & \coloneqq\bigvee\left\{ P\left(x\right)|x\in\mathcal{U}\right\} ,\\
\bigwedge_{x\in\mathcal{U}}P\left(x\right) & \coloneqq\bigwedge\left\{ P\left(x\right)|x\in\mathcal{U}\right\} .
\end{align*}
In a full fledged first order predicate logic more complicated combinations
of the two quantifiers -- referring to different variables -- occur.
Since in classical logic such more intricate quantifications are not
used, we shall largely bypass this issue.
\begin{rem*}
To indicate how such more complicated quantifications can be introduced,
we observe that if for some fixed $x_{0}\in\mathcal{U}_{0}$ and every
$\left(x,y\right)\in\mathcal{U}=\mathcal{U}_{0}\times\mathcal{U}_{1}$
\[
P\left(x,y\right)=P\left(x_{0},y\right),
\]
then $P$ is apparently independent of $x$, which can be considered
as an `all statement with respect to $x$
\[
\bigwedge_{x\in\mathcal{U}_{0}}P\left(x,y\right).
\]
 Considering $P\left(x_{0},y\right)$ as the predicate
\[
\left(x,y\right)\mapsto P\left(x_{0},y\right),
\]
we see that $\bigwedge_{x\in\mathcal{U}_{0}}P\left(x,y\right)$ is
indeed just another predicate on $\mathcal{U}$. We can now use an
existence quantifier for $y$ to obtain a proposition of the form
\[
\bigvee_{y\in\mathcal{U}_{1}}\bigwedge_{x\in\mathcal{U}_{0}}P\left(x,y\right),
\]
which coincides with $\bigvee_{y\in\mathcal{U}_{1}}P\left(x_{0},y\right)$.
This may suffice as an outline of the procedure of building more complex
quantifications. We shall not try to give a systematic construction
here, since we will not make use of such more complex cases.
\end{rem*}

There is an important subtlety in connection with the use of $\nLeftrightarrow$
hiding here. If we define
\[
\left(P\nLeftrightarrow Q\right)\coloneqq\left(P\Leftrightarrow Q\right)^{\prime},
\]
which seems natural to do, we have that $P+Q$ is perfectly well-defined
but, despite $\nLeftrightarrow$ being equal to $+$ in the propositional
logic setting, the relation $\nLeftrightarrow$ for predicates is
\emph{not }the lifted version of $+$ in predicate logic. Indeed,
\[
\left(P\nLeftrightarrow Q\right)=\left(\bigvee_{x\in\mathcal{U}}\left(P\left(x\right)\nLeftrightarrow Q\left(x\right)\right)\right),
\]
whereas
\begin{align*}
\left(P+Q\right)\left(x\right) & =\left(P\left(x\right)+Q\left(x\right)\right),\\
 & =\left(P\left(x\right)\nLeftrightarrow Q\left(x\right)\right),\:\text{ for every }x\in\mathcal{U}.
\end{align*}
This is the reason we did not include $\nLeftrightarrow$ in the list
of operations to be lifted to the predicate logic/2-algebra level.

\subsection{The Entangled Trinity}

With now having worked our way up to see the full beauty of the entangled
trinity, which may have initially a bit of a flavour of the Babylonian
confusion of tongues, but will show to be rather useful in approaching
Lewis Carroll's puzzles. We shall now try to reap some of the benefits
of our setup.

Note that there are many alternative notations around, which are,
however, more prone to produce a `mangled' trinity. We have deliberately
chosen to focus on those commonplace notations, which clearly expose
the similarities of the three facets of what we called the entangled
trinity. The following table \ref{tab:trinity} provides concise overview of the entangled
trinity as developed in the above.

\begin{table}[htbp]
\caption{The entangled trinity\label{tab:trinity}}
{\small{}%
\begin{turn}{90}
\begin{tabular}{c}
\begin{turn}{90}
{\small{}%
\begin{turn}{90}
\begin{tabular}{|c|c|c|}
\hline 
{\small 2-algebra (Boolean algebra)} & {\small set theory} & {\small predicate logic}\tabularnewline
\hline 
{\small$0$ (zero, (even))} & {\small$0$ (empty set)} & {\small$0$ (false, void)}\tabularnewline
\hline 
{\small$1$ (one, (odd))} & {\small$1$ (all of the chosen universe $\mathcal{U}$)} & {\small$1$ (true, )}\tabularnewline
\hline 
{\small$+$ (plus)} & {\small$\vartriangle$ (symmetric difference)} & {\small$\underline{\vee}$ (exclusive or)}\tabularnewline
\hline 
{\small ~'~ (successor)} & {\small ' (complement)} & {\small ~'~ , $\neg$ (negation; not)}\tabularnewline
\hline 
{\small$\,\cdot\,$,~$\min$,~$\wedge$ (times, minimum, meet)} & {\small$\cap$ (intersection)} & {\small$\wedge$ (and, L.C. uses multiplication)}\tabularnewline
\hline 
{\small$\max$,~$\vee$ (maximum, join)} & {\small$\cup$ (union)} & {\small$\vee$ (inclusive or, L.C. uses $\dagger$)}\tabularnewline
\hline 
{\small$\leq$ (less or equal to)} & {\small$\subseteq$ (sub, subset (or equal to))} & {\small$\Rightarrow$ (implies, if $\ldots$ then $\ldots$)}\tabularnewline
\hline 
{\small$=$ (equal)} & {\small$=$ (equal)} & {\small$\Leftrightarrow$ (equivalent)}\tabularnewline
\hline 
{\small$\neq$ (not equal)} & {\small$\neq$ (not equal, different)} & {\small$\nLeftrightarrow$, $\neq$ (not equivalent, not equal, different)}\tabularnewline
\hline 
{\small$\bigwedge$ (big minimum, minimum, product)} & {\small$\bigcap$ (big intersection, intersection)} & {\small$\bigwedge$ (all quantor; for every $\ldots$ we have that
$\ldots$ (is true))}\tabularnewline
\hline 
{\small$\bigvee$ (big maximum, maximum)} & {\small$\bigcup$ (big union, union)} & {\small$\bigvee$ (existence quantor; there is $\ldots$ such that
$\ldots$ (is true))}\tabularnewline
\hline 
\end{tabular}
\end{turn}}
\end{turn}\tabularnewline
\end{tabular}
\end{turn}}
\end{table}

With the induced identifications, we could indeed mix all notations
to express everything stated with the intended flavour of expression.
For example, we could write
\[
\left(A\leq B\leq C\right)\implies\left(A\leq C\right)
\]
to express the transitivity of $\leq$ as a logical implication. Or
equivalently, we could focus on logical implication altogether
\begin{equation}
\left(A\implies B\implies C\right)\implies\left(A\implies C\right),\label{eq:arrow-chain}
\end{equation}
which is, as we now know is the same statement. In a pure set perspective
this corresponds to the Venn diagram
\begin{center}
\includegraphics{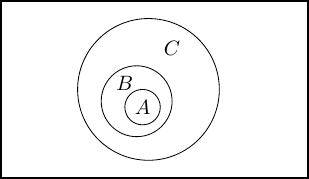}
\par\end{center}

Again, the same topic could be dressed up in the implication obtaining
the reduced equation 
\[
AC^{\prime}=0
\]
from the system of equations
\begin{align*}
AB^{\prime} & =0,\\
BC^{\prime} & =0,
\end{align*}
which is Lewis Carroll's perspective on this particular situation
(he writes $AC_{0}^{\prime},AB_{0}^{\prime},BC_{0}^{\prime},$ for
$AC^{\prime}=0,\,AB^{\prime}=0,BC^{\prime}=0$). In other terms,
\[
\left(\left(AB^{\prime}=0\right)\wedge\left(BC^{\prime}=0\right)\right)\implies\left(AC^{\prime}=0\right)
\]
is an alternative description of \eqref{eq:arrow-chain}. Since, here
something being `false' means the same as the opposite being `true',
we have
\[
\left(\left(AB^{\prime}\right)^{\prime}=1\wedge\left(BC^{\prime}\right)^{\prime}=1\right)\implies\left(\left(AC^{\prime}\right)^{\prime}=1\right)
\]

This in turn could be transposed as the equivalent implication
\[
\left(AC^{\prime}=1\right)\Longrightarrow\left(\left(AB^{\prime}=1\right)\vee\left(BC^{\prime}=1\right)\right),
\]
or in pure set language as
\[
\left(A\cap C^{\prime}\right)\subseteq\left(\left(A\cap B^{\prime}\right)\cup\left(B\cap C^{\prime}\right)\right)
\]
to possibly invoke the intuitive insight of a Venn diagram. If so
desired, we may also look at this as the 2-algebra estimate
\[
\left(AC^{\prime}\right)\leq\max\left\{ AB^{\prime},BC^{\prime}\right\} ,
\]
In any case, it suffices to consider one of the versions and obtain
proof for all of them. In principle, we could now arbitrarily mix
all the established notations without causing the results to be different.
However, for sake of clarity we will use these notations intentionally
only if we wish to convey a certain point of view by invoking the
natural language associated with them. 

\subsection{Classical logic as propositional logic on predicates}

Classical logic deals largely only with propositional logic over predicates.
Propositional logic is embedded into predicate logic via propositional
forms, which are defined on a universe $\tilde{\mathcal{U}}=2^{K}$,
$K\in\mathbb{N}$, and predicates constructed with $\vee,\wedge$,
and negation `$\:\cdots^{\prime}\:$' only. 

Substituting $K$ predicates $P_{0},\ldots,P_{K-1}$ on a universe
$\mathcal{U}$ for the $K$ variables in the propositional forms produce
predicates over $\mathcal{U}$. In this sense we obtain propositional
logic statements about predicates. More formally this is with $p$
a propositional form (over $\tilde{\mathcal{U}}$) and
\[
\mathcal{P}\coloneqq\left(P_{0},\ldots,P_{k-1}\right)\in\tilde{\mathcal{U}}^{\mathcal{U}},
\]
then
\[
P\left(x\right)\coloneqq p\left(\mathcal{P}\left(x\right)\right),\:x\in\mathcal{U}
\]
defines a predicate over $\mathcal{U}$. In particular, if the propositional
form $p$ is universally true, i.e.
\[
p=1
\]
then also 
\[
P=1
\]
regardless of what $\mathcal{P}$ actually is. This is the place where
fantasy can really develop freely yielding the type of logical puzzles
Lewis Carroll is famous for. The point is to make the resulting predicates
and the related assumptions about them plausible. Of course, nothing
prevents us from being totally unreasonable and choose for example
$\mathcal{P}=\left(P_{0},P_{1},P_{2}\right)$ with $P_{0}="chicken"$,
$P_{1}="cats"$ , $P_{2}="able\:to\:understand\:French"$, and assume
\[
\left(P_{1}\subseteq P_{2}\right)\wedge\left(0\neq P_{0}P_{1}\right)
\]
then we can safely conclude that 
\[
0\neq P_{0}P_{2},
\]
i.e. that `some chicken understand French'.

Replacing $0$ and $1$ by gap and stroke, respectively, we have indeed
that $\mathcal{P}$ is actually merely assigning a `bar code' to
the elements of $\mathcal{U}$, which are largely arbitrary and can and usually will be
 fictitious,  even -- on occasion -- nonsensical. 

On the other hand, we recall that the propositional form $x\neq0$
means $x=1$, but substituting a predicate $P$ for $x$ we get that
\[
P\neq0
\]
which is in general not the same as $P=1$ (indeed only for $P=1$
we have the equivalence\footnote{It is
\[
\left(P\neq0\right)=\left(P=0\right)^{\prime}=\left(P^{\prime}=1\right)^{\prime}\neq P.
\]
}). So, this type of statement is \emph{not} covered in this propositional
logic construction. Replacing $P$ by $P^{\prime}$, we see that we
also need to include statements like $P^{\prime}\neq0$, which is
the same as 
\[
P\neq1.
\]
 Restating this differently, we have in contrast to propositional
logic
\begin{align*}
\left(P=1\right) & \neq P.\\
\left(P=0\right) & \neq P^{\prime}.
\end{align*}
 After all
\[
P=0
\]
means $P\left(x\right)=0$ for every choice of $x\in\mathcal{U}$.
This is an `all' statement, in accordance with the lifting process,
which is also done via 'all' (for all arguments) statements. In contrast,
\[
P\neq0
\]
means $P\left(x\right)\neq0$ or $P\left(x\right)=1$ for \emph{some}
choice of $x\in\mathcal{U}$, which is an existence statement.To
avoid trouble, we shall separate the propositional-logic-over-predicates
part and the proper predicate logic statements of the form $P\neq0$
or $P\neq1$, which we may encounter in the context of classical logic. 

\subsection{Classical Aristotelian Logic: Syllogisms (Lewis Carroll's approach)}

Surely, it is already much too late to apologise to scholars, historians
and logic experts for the naive view of some mathematicians on what
we dare to label Aristotelian Logic. We are speculating, however,
that lack of sophistication in the current presentation might add
to easing understanding of major aspects of this topic. 

Aristotelian Logic deals with natural language propositions $P(A,B)$
of the following types: \textbf{none} of $A$ is in $B$ (\textbf{e}),
not none i.e. \textbf{some} of $A$ are $B$ (\textbf{i}),\textbf{ }
\begin{center}
\includegraphics{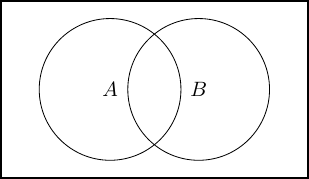}
\par\end{center}

\textbf{all} of $A$ are $B$ (\textbf{a}),
\begin{center}
\includegraphics{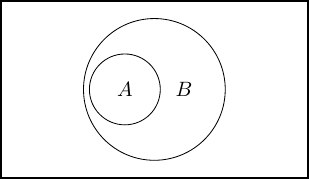}
\par\end{center}

and\textbf{ not all} of $A$ are B (\textbf{o}). These are formally
4 different propositions about $A$ and $B$ (second order predicate
logic, or just an abbreviation). The objective is to study their relationship
in order to find out which \emph{conclusion} of one of these forms
can be drawn from two other statements of this form (\emph{premises}).

Generally speaking, we have really only two types of statements (nullities)
\[
P=0
\]
and their negation (entities)
\[
P\neq0,
\]
where the special case $P=AB^{\prime}$ yields $\left(\boldsymbol{a}\right)$
and $\left(\boldsymbol{o}\right)$, see the above Venn diagrams, and
exchanging in this $B^{\prime}$ by $B$ we get cases $\left(\boldsymbol{e}\right)$
and $\left(\boldsymbol{i}\right)$.

The four basic forms can thus be formalised as
\begin{align*}
\left(A\cap B^{\prime}=0\right)
    \Leftrightarrow\left(AB^{\prime}=0\right) 
  & \Leftrightarrow\left(A\leq B\right)
    \Leftrightarrow\left(A\subseteq B\right)
    \Leftrightarrow\left(A\Rightarrow B\right)\\
  & \Leftrightarrow"\text{all \ensuremath{A} are \ensuremath{B}}",\\
\left(A\cap B=0\right)
    \Leftrightarrow\left(AB=0\right) 
  & \Leftrightarrow\left(A\leq B^{\prime}\right)
    \Leftrightarrow\left(A\subseteq B^{\prime}\right)
    \Leftrightarrow\left(A\Rightarrow B^{\prime}\right)\\
  & \Leftrightarrow"\text{none of \ensuremath{A} are \ensuremath{B}}",\\
\left(A\cap B\neq0\right)
    \Leftrightarrow\left(AB\neq0\right) 
  & \Leftrightarrow"\text{some \ensuremath{A} are \ensuremath{B}}",\\
\left(A\cap B^{\prime}\neq0\right)
    \Leftrightarrow\left(AB^{\prime}\neq0\right) 
  & \Leftrightarrow"\text{some \ensuremath{A} are not \ensuremath{B}}".
\end{align*}
Unfortunately, the tradition breaks with this neat formal symmetry
by requiring that in the first instance $A$ must not be empty. This
modification allows for the apparently desired implication
\[
"\text{all \ensuremath{A} are \ensuremath{B}}"\Rightarrow"\text{some \ensuremath{A} are \ensuremath{B}}",
\]
or more formally
\[
\left(0\neq A\right)\wedge\left(AB^{\prime}=0\right)\Rightarrow\left(AB\neq0\right),
\]
whereas
\[
\left(AB^{\prime}=0\right)\Rightarrow\left(AB\neq0\right)
\]
would be wrong for $A=0$. So, this would lead us to the modified
transcription:
\begin{align*}
\left(0\neq A\right)\wedge\left(AB^{\prime}=0\right) & \Leftrightarrow"\text{all \ensuremath{A} are \ensuremath{B}}",\\
\left(AB=0\right) & \Leftrightarrow"\text{none of \ensuremath{A} are \ensuremath{B}}",\\
\left(AB\neq0\right) & \Leftrightarrow"\text{some \ensuremath{A} are \ensuremath{B}}",\\
\left(AB^{\prime}\neq0\right) & \Leftrightarrow"\text{some \ensuremath{A} are not \ensuremath{B}}".
\end{align*}
Lewis Carroll follows this particular transcription leading him to
rephrase a simple
\[
AB^{\prime}=0
\]
(without requiring that $A$ is non-empty) as $"\text{none of \ensuremath{A} are not \ensuremath{B}}"$.
Unfortunately, this is clearly debatable, since 
\[
B=\left(B^{\prime}\right)^{\prime}
\]
shows that every $B$ can be replaced by a negation of $B^{\prime}$,
which would make it more reasonable to also have 
\begin{align*}
\left(0\neq A\right)\wedge\left(AB=0\right) & \Leftrightarrow"\text{none of \ensuremath{A} are \ensuremath{B}}".
\end{align*}
This somewhat subtle (not to say ``scholastic'') distinction is
a recurring topic of discussions with his colleagues, who tended to
ignore his choice of convention. 

In general, it is indeed false to think of $A$ as `positive' and
$A^{\prime}$ as `negative'. For example, if we wanted to label
$A$ as `being ill' and $A^{\prime}$ as `being healthy', then
$A$, i.e. `being ill', is the `positive' statement and $A^{\prime}$
its negation. But we could likewise have decided that $B$ is `being
healthy' and $B^{\prime}$ is `being ill', now making `being healthy'
the `positive' label. 

We shall bypass this questionable issue altogether by keeping the
most symmetric version above and merely recording conditions like
$0\neq A$ as a side note and then only later, in a post-processing
step, recall and utilise these types of extra assumptions, if it is
so desired. This effectively separates the propositional logic reasoning
from the genuinely predicate logic aspect in classical logic as previously
announced. 

In this context, let us briefly discuss a classical logic issue, the
so-called `problem of existence':
\begin{center}
\begin{minipage}[t]{0.87\columnwidth}%
\textbf{\emph{The problem of existence.}}

The desire to have, that if something is true for all, then it should
also be true for any particular choice is well-founded, since of course
\begin{equation}
\left(P=1\right)\Rightarrow\left(P\not=0\right).\label{eq:all-some-1}
\end{equation}
Unfortunately, the natural language is not in all cases well-suited
to distinguish between this and other, false, interpretations. Let
us consider $A\leq B$, which is the same as $AB^{\prime}=0$ or $A=AB$.
So taking $P$ to be $A=AB$ we get with \eqref{eq:all-some-1}
\[
\left(AB=A\right)=\left(\left(AB=A\right)=1\right)\Rightarrow\left(\left(AB=A\right)\not=0\right),
\]
which is certainly true. However, the natural language appears to
be not easily capable of expressing the subtle bracketing in $\left(AB=A\right)\not=0$
leading to 
\begin{equation}
AB=A\not=0,\label{eq:existence0}
\end{equation}
which is being commonly interpreted as
\[
\left(AB=A\right)\wedge\left(A\not=0\right)
\]
or
\[
\left(AB=A\right)\left(A\not=0\right).
\]
The part $A\not=0$, however, was never part of the assumption. Now,
dropping the middle term in \eqref{eq:existence0} yields
\[
AB\not=0,
\]
which is what is taken as saying `some $A$ are $B$'. This is in
a nutshell the source of the classical issue of all statements falsely
leading to an existence statement.

Matters are actually even worse. Since we have also to deal with negation,
an assumption of the form $A\not=0$ yields that for consistency we
also would need to discuss the analogous issue with $A\not=1$. The
latter is -- as far as we are aware -- not even a topic in Aristotelian
logic.%
\end{minipage}
\par\end{center}

Now, if we, in the light of this, drop the requirement of sets being
non-empty, the four so-called \emph{perfect} syllogisms are
\begin{eqnarray}
\left(\left(M\subseteq B\right)\,\wedge\:\left(A\subseteq M\right)\right) & \Rightarrow & \left(A\subseteq B\right)\,,\label{eq:aaa-2}\\
\left(\left(M\subseteq B\right)\,\wedge\:\left(0\not=AM\right)\right) & \Rightarrow & \left(0\not=AB\right)\,,\label{eq:aii-2}\\
\left(\left(M\subseteq B^{\prime}\right)\,\wedge\:\left(A\subseteq M\right)\right) & \Rightarrow & \left(A\subseteq B^{\prime}\right)\,,\label{eq:eae-2}\\
\left(\left(M\subseteq B^{\prime}\right)\,\wedge\,\left(0\not=AM\right)\right) & \Rightarrow & \left(0\not=AB^{\prime}\right)\,,\label{eq:eio-2}
\end{eqnarray}
which we shall rather use in the equivalent form
\begin{eqnarray}
\left(\left(MB^{\prime}=0\right)\,\wedge\:\left(AM^{\prime}=0\right)\right) & \Rightarrow & \left(AB^{\prime}=0\right)\,,\label{eq:aaa-2-1}\\
\left(\left(MB^{\prime}=0\right)\,\wedge\:\left(0\not=AM\right)\right) & \Rightarrow & \left(0\not=AB\right)\,,\label{eq:aii-2-1}\\
\left(\left(MB=0\right)\,\wedge\:\left(AM^{\prime}=0\right)\right) & \Rightarrow & \left(AB=0\right)\,,\label{eq:eae-2-1}\\
\left(\left(MB=0\right)\,\wedge\,\left(0\not=AM\right)\right) & \Rightarrow & \left(0\not=B^{\prime}A\right)\,.\label{eq:eio-2-1}
\end{eqnarray}
This reveals that the classical Aristotelian logic is in set terms
all about overlap of sets, such as $0\not=AB$, 
\begin{center}
\includegraphics{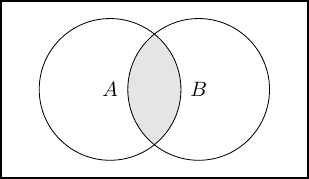}
\par\end{center}

or the lack of it, as in $0=AB$.
\begin{center}
\includegraphics{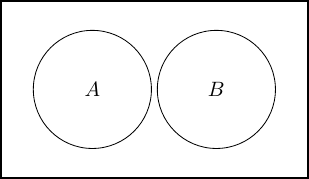}
\par\end{center}

We also note that the second set of implications are the same as the
first two with the role of $B$ replaced by $B^{\prime}$ and are
therefore redundant. This is a reflection of the redundancy in the
basic propositions $a$, $i$, $e$, $o$ (e.g. only $a$, $i$ are
really needed).

This leaves us with just two cases:
\begin{eqnarray}
\left(\left(MB^{\prime}=0\right)\,\wedge\:\left(AM^{\prime}=0\right)\right) & \Rightarrow & \left(AB^{\prime}=0\right)\,,\label{eq:aaa-2-1-1}\\
\left(\left(MB^{\prime}=0\right)\,\wedge\:\left(0\not=AM\right)\right) & \Rightarrow & \left(0\not=AB\right)\,.\label{eq:aii-2-1-1}
\end{eqnarray}

\begin{casenv}
\item \eqref{eq:aaa-2-1-1} This is the case of propositional logic about
predicates and we can apply propositional logic without even changing
notation. Indeed, we already know from the 2-arithmetic setting \eqref{eq:Ass},
\eqref{eq:Conc}, that $M$ is the eliminand, we have
\begin{equation}
\max\left\{ MB^{\prime},AM^{\prime}\right\} \geq AB^{\prime},\label{eq:2}
\end{equation}
and we cannot improve this lower estimate in terms of the retinends
$A$ and $B^{\prime}$ only. We read off that
\begin{align*}
\left(\left(MB^{\prime}=0\right)\,\wedge\:\left(AM^{\prime}=0\right)\right)=\left(\max\left\{ MB^{\prime},AM^{\prime}\right\} =0\right)\\
  \Rightarrow\left(AB^{\prime}=0\right)
\end{align*}
or
\[
\left(\left(MB^{\prime}=0\right)\,\wedge\:\left(AM^{\prime}=0\right)\right)\leq\left(AB^{\prime}=0\right).
\]
As a Venn diagram
\begin{center}
\includegraphics{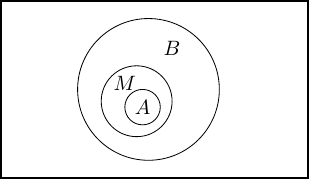}
\par\end{center}

or as a graph of inclusions
\begin{center}
\includegraphics{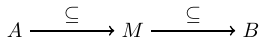}
\par\end{center}
\item \eqref{eq:aii-2-1-1} This case involves a $\neq$ statement, which
we wanted to set aside. Noting that $AM\leq A,\:AM\leq M$ we introduce
a dummy predicate $X\neq0$ with
\begin{equation}
\left(X\leq A\right)\wedge\left(X\leq M\right),\label{eq:=000020dummy}
\end{equation}
which is equivalent to $0\neq X\leq AM$. We know that such an $X$
exists from the original assumption (without the dummy predicate $X$)
and conversely, the existence of such an $X$ implies the original
assumption. In terms of nullities \eqref{eq:=000020dummy} is
\[
\left(XA^{\prime}=0\right)\wedge\left(XM^{\prime}=0\right).
\]
Thus, according to \eqref{eq:aii-2-1-1} the assumption becomes 
\[
\left(\left(XA^{\prime}=0\right)\wedge\left(\left(XM^{\prime}=0\right)\wedge\left(MB^{\prime}=0\right)\right)\right)
\]
with now $X,A,B^{\prime}$ being retinends and $M$ still being the
eliminand. The situation as a Venn diagram
\begin{center}
\includegraphics{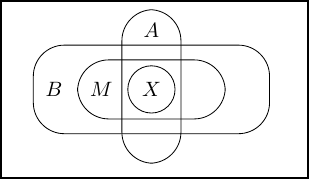}
\par\end{center}

or as a graph of inclusions
\begin{center}
\includegraphics{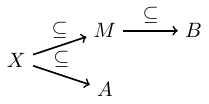}
\par\end{center}

The term $\left(XM^{\prime}=0\right)\wedge\left(MB^{\prime}=0\right)$
in the extra -- normally unnecessary -- brackets are assumptions
of type \eqref{eq:aaa-2-1-1} yielding as in Case 1
\[
\max\left\{ MB^{\prime},XM^{\prime}\right\} \geq\left(XB^{\prime}\right).
\]
 The remaining term $XA^{\prime}=0$ is left unchanged, since it only
contains retinends. So, we have altogether
\begin{equation}
\max\left\{ XA^{\prime},MB^{\prime},XM^{\prime}\right\} \geq\max\left\{ XA^{\prime},XB^{\prime}\right\} \label{eq:3}
\end{equation}
and so 
\begin{align*}
\left(\left(XA^{\prime}=0\right)\,\wedge\:\left(MB^{\prime}=0\right)\,\wedge\:\left(XM^{\prime}=0\right)\right)\\
  \leq\left(\max\left\{ XA^{\prime},XB^{\prime}\right\} =0\right).  
\end{align*}
Now, finally we want to eliminate the dummy predicate $X$. Noting
that
\[
\left(\max\left\{ XA^{\prime},XB^{\prime}\right\} =0\right)=X\:\max\left\{ A^{\prime},B^{\prime}\right\} =0
\]
\[
\left(\left(XA^{\prime}=0\right)\,\wedge\:\left(MB^{\prime}=0\right)\,\wedge\:\left(XM^{\prime}=0\right)\right)\leq\left(X\max\left\{ A^{\prime},B^{\prime}\right\} =0\right).
\]
 We find
\begin{align*}
\left(X\max\left\{ A^{\prime},B^{\prime}\right\} =0\right)=\left(X\left(\min\left\{ A,B\right\} \right)^{\prime}=0\right)\\
  =\left(X\left(AB\right)^{\prime}=0\right)
  =\left(X\leq\left(AB\right)\right).
\end{align*}
At this this point we recall $X\neq0$ and note 
\[
0\neq X\leq\left(AB\right),
\]
which implies
\[
0\neq AB.
\]
Thus, we have shown \eqref{eq:aii-2-1-1}. We have been rather detailed
in this particular derivation, since we take this case as an illustration
how to handle $\neq$ statements in general, should they occur.
\end{casenv}
As a first example we consider a just about `proverbial application'
of classical Aristotelian logic, which is the syllogism:
\begin{enumerate}
\item Socrates is human.
\item Humans are mortal.
\end{enumerate}
The conclusion is 

\begin{enumerate}

\setcounter{enumi}{2}

\item Socrates is mortal.

\end{enumerate}

This syllogism has its issues. So to begin with, to make this meaningful,
there appears to be a universe of `beings' needed, where `immortals'
are existing\footnote{Since `immortals' do not exist, to say `somebody is mortal' is
actually an empty statement. On the other hand, in natural language
sentences such as `we are all mortal', which are equally void, are
commonly used. In any case, the argument still works.}. Lets say $\mathcal{U}$ is all `humans and gods' then. 

There has also been a debate that this is not a proper syllogism at
all, since Sokrates is the name of an individual not a predicate.
As it turns out this is not a logic problem but just a natural language
issue, since in natural language it is not commonplace to distinguish
between Sokrates as an element of $\mathcal{U}$ and the singleton
$\left\{ Sokrates\right\} $ as a subset of $\mathcal{U}$. Of course
in the translation to formal logic terms we would choose the latter
transcription. So, with $H$ denoting humans, $M$ for mortals, $S$
for $\left\{ Sokrates\right\} $ we transcribe the syllogism to
\begin{enumerate}
\item $S\subseteq H$,
\item $H\subseteq M$,
\item $S\subseteq M$, 
\end{enumerate}
which is of type \eqref{eq:aaa-2-1-1}. Since $S\subseteq H$ and
$S\subseteq M$ is, due to $S$ having only one element, is actually
equivalent to $Sokrates\in H$ or $Sokrates\in M$, respectively,
the natural language translation is perfectly recovered. 

\subsubsection*{The Jack Sprat Problem}
\begin{center}
\textbf{\emph{\textquotedblleft Jack Sprat could eat no fat: His wife
could eat no lean: And so, between them both, They licked the platter
clean.\textquotedblright{}}}
\par\end{center}

So, the poem starts from a situation, where Jack and his wife have
been served with a platter of fat and lean meat (as proposed by Lewis
Carroll, see Bartley, Book XIII, Chapter X, \cite{book834284}). The
problem that Lewis Carroll posed is to set the poem in a suitable
frame so that it can be expressed in admissible syllogistic terms. 

To express the intended conclusion ''And so, between them both, They
licked the platter clean'' formally leads first to
\[
\mathcal{U}=1=J^{\prime}\cup W^{\prime},
\]
where 
\[
\mathcal{U}\text{="meat on the platter served",}
\]
 
\[
\text{\ensuremath{J=}"meat from the platter served, which Jack will not eat"}
\]

and 

\[
W=\text{"meat from the platter served, which Jack's wife will not eat"}.
\]
There is an implicit assumption here that ``cannot'' means ``will
not'', so that $J^{\prime}$ and $W^{\prime}$ conveniently encodes
\textquotedbl meat from the platter served, which Jack will eat\textquotedbl{}
and \textquotedbl meat from the platter served, which Jack's wife
will eat\textquotedbl , respectively. These are the kind of liberties
we find in Lewis Carroll's examples all the time to allow for a more
varied, concise and entertaining wording of his riddles.

To continue with the task, recall that we can, by taking complements,
reformulate the desired conclusion with the abbreviations provided
as
\[
0=JW,
\]
which in turn is equivalent to
\begin{equation}
J\subseteq W^{\prime}.\label{eq:jack-conclusion}
\end{equation}
The remaining task is to obtain \eqref{eq:jack-conclusion} as the
conclusion from a syllogism. In other words, we aim to find an $M$
such that
\[
J\subseteq M\subseteq W^{\prime},
\]
where $J\subseteq M$ and $M\subseteq W^{\prime}$ is supposed to
represent the premises. In the light of the verbalisation of the premises
we are led to take $M$=''fat meat on the platter served''. In summary
the poem verbalises the syllogism\footnote{Recall that
\[
\left(J\subseteq W^{\prime}\right)\Leftrightarrow\left(J^{\prime}\cup W^{\prime}=1\right).
\]
}
\[
J\subseteq M\subseteq W^{\prime}\Rightarrow J\subseteq W^{\prime}.
\]
It is of course debatable if the syllogism should not rather be
\[
0\not=J\subseteq M\subseteq W^{\prime}\not=1\Rightarrow0\not=J\subseteq W^{\prime}\not=1.
\]
In any case is this a minimal alternative solution to the one by T.
W. Settle reported in Book XIII, Chapter X, \cite{book834284}: 
\begin{itemize}
\item Premise 1: The meat on the platter that Jack will not eat is fat meat. 
\item Premise 2: Fat meat is what Jack's wife will eat. 
\item Conclusion: All meat Jack does not eat, his wife eats, which is equivalently
stated ''And so, between them both, They licked the platter clean''.
\end{itemize}

\subsection{Classical Aristotelian Logic: Sorites\! (Lewis Carroll's approach\!)}

Since\! \eqref{eq:aaa-2-1-1} leads to \eqref{eq:2}, which has an assumption
with two products of two factors, and \eqref{eq:aii-2-1-1} leads
to \eqref{eq:3}, which has an assumption with three products of two
factors, there is an obvious generalisation to assumptions that are a maximum of an arbitrary number of products of two factors. This is the basis of a large class
of so-called \emph{sorites} (``heap of syllogisms'', syllogismos
synthetos). There are simpler cases, which are just a long chain of
inclusions
\begin{center}
\includegraphics{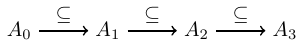}
\par\end{center}

or as a Venn diagram
\begin{center}
\includegraphics{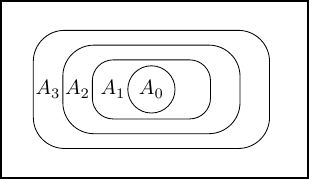}
\par\end{center}

or more complicated situations. To illustrate the latter, let us just
look at an arbitrary graph of inclusions
\begin{center}
\includegraphics{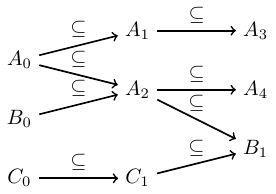}
\par\end{center}

which generalised the simpler chain case as well as the graphs encountered
in the context of syllogisms. We shall only look more thoroughly at
some examples to illustrate the method.

\subsubsection{Babies and Crocodiles}

For starters, let us reconsider the example about the babies and crocodiles
briefly mentioned in the introductory section \ref{subsec:The-idea}. 

(1) Babies are illogical;

(2) Nobody is despised who can manage\emph{ }a crocodile; 

(3) Illogical persons are despised.

What should we use as a universe $\mathcal{U}$? Well, let us pick
any universe containing people (which should include babies), maybe
$\mathcal{U}$ could be just `people'. There are predicates
or subsets to consider: $B$ for babies, $D$ for despised elements
of $\mathcal{U}$, $I$ for illogical elements of $\mathcal{U}$,
$C$ for elements of $\mathcal{U}$ able to manage crocodiles. Then
we translate (1)-(3) into inclusion statements:
\begin{enumerate}
\item $B\subseteq I$,
\item $C\subseteq D^{\prime}$ (or $D\subseteq C^{\prime})$,
\item $I\subseteq D$.
\end{enumerate}
Indeed, combining them we see we have a single chain
\[
B\subseteq I\subseteq D\subseteq C^{\prime}.
\]
Writing this in terms of nullities, this is
\[
BI^{\prime}=ID^{\prime}=DC=0.
\]
We see that $I,D$ are eliminands and $B,C$ are retinends. So, we
suspect
\begin{equation}
BC=0\label{eq:baby-conc}
\end{equation}
 to be the conclusion. Indeed, we find
\begin{equation}
\max\left\{ BI^{\prime},ID^{\prime},DC\right\} \geq\left(BC\right)\label{eq:3x}
\end{equation}
with $BC$ being the largest lower bound. We just assume $BC=1$,
which -- since we work in 2-arithmetic -- yields $B=C=1$ and find
\[
\max\left\{ I^{\prime},ID^{\prime},D\right\} =\begin{cases}
\max\left\{ I^{\prime},I\right\} =1 & \text{for }D=0,\\
\max\left\{ I^{\prime},0,1\right\} =1 & \text{for }D=1,
\end{cases}
\]
confirming the inequality. The same way we see, that having just either
$B=1$ or $C=1$ is not enough to force the inequality to be always
true. So, $BC$ really is the largest lower bound and so \eqref{eq:baby-conc}
is the conclusion. This is still an example with three products of
two factors.

But we see that if we know that the assumption has a certain number
of products of two factors, we also find that the largest lower bound
can also only have a number of products of two factors. Thus, we only
need to consider all products of two retinends and determine if they
are lower bounds. Then the conclusion is that all these lower bounds
together set equal to $0$ is the conclusion. We shall employ this
observation to solve the following larger puzzles.

\subsubsection{Animals}

Charles Ludwidge Dodgson has raised the question how to solve larger
sorites in his book on logic. Here is a more involved example to illustrate
the outlined procedure applied to sorites. 

The assumptions are
\begin{enumerate}
\item Alligators are reptiles.
\item Some mammals are carnivorous.
\item No reptile is a mammal.
\item Jumbo is an elephant.
\item No mammals are alligators.
\item No elephant is an alligator.
\end{enumerate}
The task is to find the conclusion.

Taking `animals' as the universe $\mathcal{U}$ and $A$ as the
alligators, $M$ a mammals, $C$ as carnivores, $R$ as reptiles,
$E$ as elephants, $J$\footnote{Recall that there is a classical discussion about the issue if individuals
are permitted in a syllogism. In our view this issue is resolved by
pointing out that here we are not talking about a particular elephant
by the name ``Jumbo'' but of a set $J$ containing only one element
$J\coloneqq\left\{ "Jumbo"\right\} $. One must not confuse $x$ and
$\left\{ x\right\} $!} for $\left\{ "Jumbo"\right\} $, this translates first to 
\begin{multicols}{2}

\begin{enumerate}
\item $0\not=A$ and $A\subseteq R$, 
\item $0\not=CM$, 
\item $M\subseteq R^{\prime}$, 
\item $0\not=J\subseteq E$, 
\item $M\subseteq A^{\prime}$, 
\item $E\subseteq A^{\prime}$. 
\end{enumerate}
\end{multicols}
We find $A,X,J\neq0$, where $X$ is the dummy variable to resolve
the `some' statement as above. Storing this aside we are left with
the translation into nullities: 
\begin{multicols}{2}

\begin{enumerate}
\item $AR^{\prime}=0,$ 
\item $XC^{\prime}=0$, $XM^{\prime}=0$, 
\item $MR=0,$ 
\item $JE^{\prime}=0$, 
\item $MA=0,$ 
\item $EA=0.$ 
\end{enumerate}
\end{multicols}
 So this would be an assumption with 7 products with two factors.
Since we call a predicate a retinend, if it never appears both complemented and uncomplemented among the factors of the products in the assumptions, we read off that $A,X,J,C^{\prime}$ are the only retinends and so $R,M,E$ are
eliminands. Since we know that only two factors need to be tested
we end up with 6 possible largest lower bounds of $\max\left\{ AR^{\prime},XC^{\prime},XM^{\prime},MR,JE^{\prime},MA,EA\right\} $
in product form, namely 
\[
AX,AJ,AC^{\prime},XJ,XC^{\prime},JC^{\prime}.
\]

\begin{casenv}
\item $AX$\\
 Let $AX=1$. Since 
\begin{align*}
0 & =\max\left\{ R^{\prime},C^{\prime},M^{\prime},MR,JE^{\prime},M,E\right\} \\
 & =\max\left\{ C^{\prime},M^{\prime},M,JE^{\prime},M,E\right\} 
\end{align*}
has no solution, we confirm that $AX$ is a lower bound. 
\item $AJ$\\
 Let $AJ=1$. Since 
\[
\max\left\{ R^{\prime},XC^{\prime},XM^{\prime},MR,E^{\prime},M,E\right\} =0
\]
has no solution, we confirm that $AJ$ is a lower bound. 
\item $AC^{\prime}$\\
 Let $AC^{\prime}=1$. Since 
\begin{align*}
0 & =\max\left\{ R^{\prime},X,XM^{\prime},MR,JE^{\prime},M,E\right\} \\
 & =\max\left\{ J\right\} .
\end{align*}
has $J=0$ as solution, this is not a lower bound. 
\item $XJ$\\
 Let $XJ=1$. Since 
\begin{align*}
0 & =\max\left\{ AR^{\prime},C^{\prime},M^{\prime},MR,E^{\prime},MA,EA\right\} \\
 & =A
\end{align*}
has $A=0$ as solution, this is not a lower bound. 
\item $XC^{\prime}$\\
 Let $XC^{\prime}=1$. Since 
\begin{align*}
0 & =\max\left\{ AR^{\prime},1,M^{\prime},MR,JE^{\prime},MA,EA\right\} 
\end{align*}
has no solution, we confirm that $XC^{\prime}$ is a lower bound. 
\item $JC^{\prime}$\\
 Let $JC^{\prime}=1$. Since 
\begin{align*}
0 & =\max\left\{ AR^{\prime},X,XM^{\prime},MR,E^{\prime},MA,EA\right\} \\
 & =MR
\end{align*}
has $M=0$ or $R=0$ as solution, this is not a lower bound. 
\end{casenv}
Thus, we have that 
\[
\max\left\{ AX,XC^{\prime},AJ\right\} =0
\]
is the conclusion. It remains to eliminate the dummy variable $X$
leading us to 
\[
\left(0\neq\left(A^{\prime}C\right)\right)\wedge J\subseteq A^{\prime}.
\]

This can be verbalised as ``Jumbo is not an alligator, and there are carnivores that are not alligators.''

\subsubsection{Cats, Mice, and Kangaroos}

The following is taken from \cite{book834284}, p. 174:

(1) The only animals in this house are cats;

(2) Every animal is suitable for a pet, that loves to gaze at the
moon;

(3) When I detest an animal, I avoid it;

(4) No animals are carnivorous, unless they prowl at night;

(5) No cat fails to kill mice;

(6) No animals ever take to me, except what are in this house;

(7) Kangaroos are not suitable for pets;

(8) None but carnivora kill mice;

(9) I detest animals that do not take to me;

(10) Animals, that prowl at night, always love to gaze at the moon.

$\mathcal{U}$= \textquotedbl animals\textquotedbl ; $a$ = avoided
by me; $b$ = carnivora; $c$ = cats; $d$ = detested by me; $e$=in
this house; $h$ = kangaroos; $k$ = killing mice; $l$ = loving to
gaze at the moon; $m$ = prowling at night; $n$ = suitable for pets;
$r$ = taking to me. The following Venn diagram describes the situation
(dropping $ar^{\prime}\subseteq d\subseteq a$ for simplification)
\begin{center}
\includegraphics{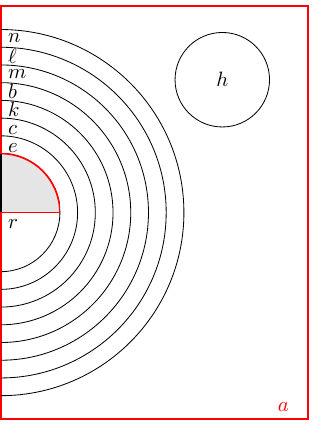}
\par\end{center}

The task is as usual to find the conclusion. The first step in this
setup is as usual a conversion into formulas.

\begin{multicols}{2}

\begin{enumerate}
\item $0\not=e\subseteq c$, 
\item $0\not=l\subseteq n$, 
\item $d\subseteq a$, 
\item $b\subseteq m$ , 
\item $c\subseteq k$, 
\item $r\subseteq e$, 
\item $0\not=h\subseteq n^{\prime}$, 
\item $k\subseteq b$ , 
\item $r^{\prime}\subseteq d$, 
\item $0\not=m\subseteq l$, 
\end{enumerate}
\end{multicols}

or, recording $0\not=e,0\not=l,0\not=h,0\not=m,$ for later use we
are left with the nulltities

\begin{multicols}{2}

\begin{enumerate}
\item $ec^{\prime}=0$, 
\item $ln^{\prime}=0$, 
\item $da^{\prime}=0$, 
\item $bm^{\prime}=0$ , 
\item $ck^{\prime}=0$, 
\item $re^{\prime}=0$, 
\item $hn=0$, 
\item $kb^{\prime}=0$ , 
\item $r^{\prime}d^{\prime}=0$, 
\item $ml^{\prime}=0$. 
\end{enumerate}
\end{multicols}

Since only $a^{\prime},h$, are retinends, we have only one candidate
for being the largest lower bound of 
\[
\max\left\{ ec^{\prime},ln^{\prime},da^{\prime},bm^{\prime},ck^{\prime},re^{\prime},hn,kb^{\prime},r^{\prime}d^{\prime},ml^{\prime}\right\} 
\]
namely $a^{\prime}h$. Indeed, assuming $a^{\prime}h=1$ we have to
inspect 
\begin{align*}
0 & =\max\left\{ ec^{\prime},ln^{\prime},d,bm^{\prime},ck^{\prime},re^{\prime},n,kb^{\prime},r^{\prime}d^{\prime},ml^{\prime}\right\} \\
 & =\max\left\{ b,k^{\prime},b^{\prime}\right\} ,
\end{align*}
finding that there is no solution confirming, that $a^{\prime}h$
is the largest lower bound. This yields the conclusion 
\[
a^{\prime}h=0,
\]
which can be verbalised as ''I avoid kangaroos''.

If products of more than two factors occur in the assumption our shortcut
cannot be used anymore and we therefore shall now lay out the general
procedure. With the general procedure established we shall then inspect
some more of Lewis Carroll's puzzles in the framework developed.

\section{The General Setting}

\subsection{Choosing a Suitable Universe}

We have seen several cases of universes $\mathcal{U}$ chosen. In
the classical context these are always finite. Universes are, at this informal stage, pre-set-theoretic and pre-logical objects: simply collections of `things' that our predicates are meant to talk about. They belong to a naive
set concept rooted in the realm of natural language rather than mathematical
theory. There may even be `things' in the universe that do not exist in any real sense, but are merely placeholders or names we want to apply predicates to.

The earlier example `John Sprat' already shows the importance of choosing
the right universe. In some examples there may actually be the need
to choose a sufficiently small universe. There is of course always
the smallest universe $\mathcal{U}$ consisting of just the union
of all verbal predicates involved (interpreted as sets) as a fallback.
In any case, due to the unclear terminology of natural languages there
is opportunity of discussion, which universe is the right one to choose.

An example, where the smallness of $\mathcal{U}$ is an issue, would
be if we find out that to make things work we need `not daughters'
to mean `sons' as in the following example taken from ``Game of Logic'',
Chapter I, \textbf{79}, example 29, see e.g. \cite{Carroll1958}:

(1) All my sons are slim;

(2) No child of mine is healthy who takes no exercise;

(3) All gluttons, who are children of mine, are fat;

(4) No daughter of mine takes any exercise.

The suggested abbreviations are $a$ = ``fat''; $b$ = ``gluttons'';
$c$ = ``healthy''; $d$ = ``sons''; $e$ = ``taking exercise''.
Here the universe should be taken as $\mathcal{U}$= ``my children,
which are slim or fat'' to ensure that $a^{\prime}=\text{"slim"}$
$d^{\prime}=\text{"daughters"}$. Then (1)-(4) translates to 
\begin{multicols}{2}

\begin{enumerate}
\item $d\subseteq a^{\prime}$, 
\item $c\subseteq e$, 
\item $b\subseteq a$, 
\item $d^{\prime}\subseteq e^{\prime}.$ 
\end{enumerate}
\end{multicols}

In terms of nullities the assumption is 
\[
\max\left\{ ad,ce^{\prime},ba^{\prime},ed^{\prime}\right\} =0.
\]
Using the shortcut of the previous section, the suspected largest
lower bound of $\max\{ ad,\\ ce^{\prime}, ba^{\prime}, ed^{\prime}\}$
is $bc$. Indeed, letting $bc=1$ we inspect 
\begin{align*}
0 & =\max\left\{ ad,e^{\prime},a^{\prime},ed^{\prime}\right\} \\
 & =\max\left\{ d,d^{\prime}\right\} ,
\end{align*}
which has no solution, thus confirming that 
\[
bc=0
\]
is the conclusion. This may be verbalised as: ``my healthy children
are no gluttons'' or ``the gluttons among my children are not healthy''.

\subsection{Eliminands and Retinends}

We recall from earlier discussion that assumptions need to be turned
into a `normal form' as a maximum of pro\-ducts set equal to $0$, where
the products need to be tailored as close as possible to the natural
language formulation\footnote{This may be at the heart of Lewis Carroll's criticism of Boole's approach
to logic, which he criticises for being too abstract and unintuitive.
The closeness to natural language is essential in determining what
eliminands and retinends are. One could consider to construct minimal
disjunctive normal forms instead, but this would be a hard problem
and way beyond, what we want to achieve here.}. If needed we introduce dummy variables to make `$0\neq$' statements
better digestible. We noted that a logically equivalent assumption
may have a completely different set of retinends and eliminands! So,
this issue is still associated with natural language and not yet fully
in the realm of logic. In consequence, there is also plenty of opportunity
to debate, which is the right normal form to use. We recall that retinends
are all those predicates appearing as factors in the product of the
normal form of the assumption, which are not also occurring in negated
form. So, if $A$ occurs in the products but there is no $A^{\prime}$
appearing then $A$ is a retinend. Similarly, if $A^{\prime}$ occurs
in the products but there is no $A$ appearing then $A^{\prime}$
is a retinend. All others are eliminands. In other words, whenever both $A$
and $A^{\prime}$ appear as factors among the products in the normal form of the assumption, then 
$A$ is called an eliminand. After having identified
retinends and eliminands, we are, however, safely in the realm of
logic.

\subsection{Proper Complete Conclusions}

To our knowledge Lewis Carroll never explicitly defined, what a `complete
conclusion' is supposed to be. Our suggestion here is to determine
the largest lower bound of $\mathcal{A}$ (as it occurs in the assumption
$\mathcal{A}=0$) that can be obtained as the maximum of largest lower
bounds in terms of retinends only. This way we have a unique outcome.
The procedure is extrapolating the construction for syllogisms. As
a general observation we record that a conclusion is obtained by \emph{removing}
information from the assumptions in the sense that the conclusion corresponds to a weaker statement, i.e. a larger set in the inclusion order. So, we see
that what is implied cannot be more specific than the assumption.
A proper inclusion always goes along with a more vague proposition
being implied from the more specific assumption. Still, to our limited
mind the conclusion appears to convey a new meaning, since the assumption
was distracting or confusing us. 

There are some fine points to consider to post-process the unique
outcome found as the largest lower bound in terms of retinends only.
The first step is the removal of dummy variables introduced in the
process of formalisation (as illustrated above). After that we are
permitted to speak of \emph{the} complete conclusion. There is, however,
one more step to perform to come to what we would want to call the
\emph{proper} complete conclusion. It may be that the complete conclusion
contains a term, which just repeats a product occurring in $\mathcal{A}$.
Apparently, Lewis Carroll would prefer to strike out such terms, since
there is no new insight coming from them. Also, largest lower bounds
equivalent to $0$ or $\mathcal{A}$ would not be accepted as conclusions
for the same reason (according to a letter to John Cook Wilson, see
\foreignlanguage{english}{\cite{book834284}}, p. 374). The lower
bound $0$ is trivial and a lower bound logically equivalent to $\mathcal{A}$
merely repeats the assumption altogether. The latter could be hidden
in finding a best lower bound of $\mathcal{A}$ in terms of retinends
only to be $1$, which would simply say that the assumption $\mathcal{A}$
is universally true already, i.e. $\mathcal{A}=1$. 

As an example, how Lewis Carroll struggled with the term conclusion,
we consider the `Library Problem'.

\subsubsection*{The Library Problem}

The problem goes like this:

(1) All the old books are Greek;

(2) All the quartos are bound;

(3) None of the poets are old quartos.

$\mathcal{U}$=''books in the library'', $a$=''bound books'',
$b$=''Greek books'', $c$=''old books'',\! $e$=''Quartos'',
$d$=''books on poetry''. So, this translates to 
\begin{multicols}{3}

\begin{enumerate}
\item $0\neq c\subseteq b$, 
\item $0\neq e\subseteq a$, 
\item $d\subseteq e^{\prime}\cup c^{\prime}$ 
\end{enumerate}
\end{multicols}

or setting $0\neq c,e$ aside and listing the rest as nullities we
get 
\begin{multicols}{3}

\begin{enumerate}
\item $cb^{\prime}=0$, 
\item $ea^{\prime}=0$, 
\item $ced=0$. 
\end{enumerate}
\end{multicols}

Clearly there are no eliminands and so here nothing can be eliminated!
John Cook Wilson sent a proposed solution to the problem, which Lewis
Carroll found to be equivalent to the assumption and so not acceptable
as a proper conclusion. As an afterthought he noticed that indeed
there may not be an acceptable conclusion. There are also extra\footnote{There is actually a calculation given by Lewis Carroll in this letter
supposedly solving the library problem, where $a=b=c=0$ is imposed
additionally to $a^{\prime}b^{\prime}c^{\prime}d^{\prime}e^{\prime}=0$.
It is unclear what this is intended to show, but the result is that
$a^{\prime}b^{\prime}c^{\prime}d^{\prime}e^{\prime}=0$ reduces to
\begin{align*}
d^{\prime}e^{\prime} & =0.
\end{align*}
} constraints to the library problem, which are only mentioned in the
letter and missing from the recorded setup, namely that $a,b,c,d,e$
cover all books in the library, i.e.

(4) $a^{\prime}b^{\prime}c^{\prime}d^{\prime}e^{\prime}=0$.

After transforming the assumptions into the chosen normal form, only 
$c,d,e$ occur in both polarities among the products, hence they are eliminands, while 
$a^{\prime},b^{\prime}$ remain retinends. So, let $a^{\prime}b^{\prime}=1$ then $\max\left\{ a^{\prime}b^{\prime}c^{\prime}d^{\prime}e^{\prime},cb^{\prime},ea^{\prime},cde\right\} $
reduces 
\begin{align*}
0 & =\max\left\{ a^{\prime}b^{\prime}c^{\prime}d^{\prime}e^{\prime},cb^{\prime},ea^{\prime},cde\right\} ,\\
 & =\max\left\{ c^{\prime}d^{\prime}e^{\prime},c,ea^{\prime},cde\right\} ,\\
 & =\max\left\{ d^{\prime}e^{\prime},ea^{\prime}\right\} ,
\end{align*}
which clearly has a solution, showing that the only candidate available
is not a lower bound. Leaving us without a proper complete solution
also in this case. In the correspondence with John Cook Wilson there
is another constraint mentioned, which amounts to $a^{\prime}b^{\prime}c^{\prime}=1$.
Taking that as the test case for a lower bound (ignoring that $c$
is an eliminand) we would find that 
\begin{align*}
0 & =\max\left\{ a^{\prime}b^{\prime}c^{\prime}d^{\prime}e^{\prime},cb^{\prime},ea^{\prime},cde\right\} ,\\
 & =d^{\prime},
\end{align*}
has the obvious solution $d=1$. So, still no nontrivial lower bound.

\subsubsection*{Wilson's problem number 9}

In a letter Lewis Carroll refers to a problem number 9 John Cook Wilson
has sent to Lewis Carroll, which is given as a formal problem as the
assumption 
\[
\mathcal{A}=0
\]
with 
\[
\mathcal{A}\coloneqq\max\left\{ AB^{\prime}C^{\prime},BCKD^{\prime}\right\} .
\]
Clearly we may reduce the distraction in this problem by introducing
just one predicate $W\coloneqq KD^{\prime}$ to obtain 
\[
\mathcal{A}\coloneqq\max\left\{ AB^{\prime}C^{\prime},BCW\right\} .
\]
A candidate for a proper lower bound as a product of retinends only
would be $AW$. So letting $AW=1$ we test 
\begin{align*}
0 & =\max\left\{ AB^{\prime}C^{\prime},BCW\right\} ,\\
 & =\max\left\{ B^{\prime}C^{\prime},BC\right\} ,
\end{align*}
which has for example $B=1$ and $C=0$ as solution. Thus, $AW$ is
not a lower bound at all.

An example problem showing a serious issue with the classical assumption
that `all $A$ are $B$' encompasses `some $A$ are $B$' is the problem
of the grocers on bikes.

\subsubsection*{The Problem of Grocers on Bicycles}

(1) All honest industrious men are healthy,

(2) No grocers are healthy,

(3) All industrious grocers are honest,

(4) All cyclists are industrious,

(5) All unhealthy cyclists are dishonest.

The conclusion is posed as `Grocers do not ride cycle'. The question
raised is to find the superfluous premise.

Let us check the conclusion first. Taking $\mathcal{U}="people"$
as our universe and $H="honest\:people"$, $I="industrious\:people"$,
$X="healthy\:people"$, $G="grocers",$ $C="cyclists"$, we have as
a first translation 
\begin{multicols}{2}

\begin{enumerate}
\item $0\neq HI\subseteq X$, 
\item $G\subseteq X^{\prime}$, 
\item $\underline{
0\neq IG}\subseteq H$, 
\item $0\neq C\subseteq I$, 
\item $0\neq X^{\prime}C\subseteq H^{\prime}$.\\
 
\end{enumerate}
\end{multicols}

From (1)\&(2) we get 
\[
0\neq HI\subseteq X\subseteq G^{\prime},
\]
which by (3) becomes 
\begin{align}
0 & \neq IG\subseteq HI\subseteq X\subseteq G^{\prime}\subseteq\left(GI\right)^{\prime}.\label{eq:grocers-chain0}
\end{align}
This yields a contradiction! To repair this we need to drop the non-emptiness
constraint here (the offending condition underlined) and obtain $IG=0$.
Following Lewis Carroll's habit, for this we would need to reword
(3) to ``none of the industrious grocers is dishonest''. One could
of course argue that only the non-emptiness $0\neq C$ should be retained,
since all other `all' statements refer to intersections and would
mean `some honest people are industrious' and `some cyclists are not
healthy'. This would, however, be a consistency problem with the concept
that `all' is supposed to imply `some'. For simplicity we ignore the
`$0\neq$' part altogether and list the remaining nullities as 
\begin{multicols}{2}

\begin{enumerate}
\item $HIX^{\prime}=0$, 
\item $GX=0$, 
\item $IGH^{\prime}=0$, 
\item $CI^{\prime}=0$, 
\item $X^{\prime}CH=0$.\\
 
\end{enumerate}
\end{multicols}
 We detect $H,X,I$ as eliminands and $G,C$ as retinends. The candidate
for a lower bound of $\max\left\{ HIX^{\prime},GX,IGH^{\prime},CI^{\prime},X^{\prime}CH\right\} $
in terms of $G,C$ is $GC$. So let $GC=1$ then 
\begin{align*}
0 & =\max\left\{ HIX^{\prime},GX,IGH^{\prime},CI^{\prime},X^{\prime}CH\right\} ,\\
 & =\max\left\{ HIX^{\prime},X,IH^{\prime},I^{\prime},X^{\prime}H\right\} ,\\
 & =\max\left\{ 1,X,X^{\prime}\right\} ,
\end{align*}
which has no solution, confirming that $GC$ is a lower bound. Now,
$G$ alone is not a lower bound, since for $G=1$ we have that 
\begin{align*}
0 & =\max\left\{ HIX^{\prime},GX,IGH^{\prime},CI^{\prime},X^{\prime}CH\right\} ,\\
 & =\max\left\{ HI,IH^{\prime},CI^{\prime},CH\right\} ,
\end{align*}
has a solution $I=C=0$. Similarly, $C$ alone cannot be a lower bound,
since for $C=1$ we have that 
\begin{align*}
0 & =\max\left\{ HIX^{\prime},GX,IGH^{\prime},I^{\prime},X^{\prime}H\right\} ,\\
 & =\max\left\{ HX^{\prime},GX,GH^{\prime}\right\} ,
\end{align*}
has a solution $G=H=0$. So the complete conclusion is 
\[
GC=0,
\]
which can indeed be verbalised as `grocers do not ride bicycles'.
For $GC=1$ we also have 
\[
1=\max\left\{ HIX^{\prime},X,IH^{\prime},I^{\prime},X^{\prime}H\right\} .
\]
The question is which assumption is superfluous. Since the contradiction
$\max\left\{ I,I^{\prime}\right\} =0$ of $GC$ not being a lower
bound was based on $I$, we see that $H^{\prime}=0$ still lead to
\begin{align*}
1 & =\max\left\{ HIX^{\prime},X,IH^{\prime},I^{\prime},X^{\prime}H\right\} ,\\
 & =\max\left\{ IX^{\prime},X,I,I^{\prime},X^{\prime}\right\} ,
\end{align*}
which leads to suspecting $HXI^{\prime}=0$ or $X^{\prime}CH=0$ being
superfluous. So, let $GC=1$ then 
\begin{align*}
0 & =\max\left\{ HIX^{\prime},GX,IGH^{\prime},CI^{\prime}\right\} ,\\
 & =\max\left\{ H,H^{\prime}\right\} ,
\end{align*}
does not have a solution and so confirms $GC$ as a lower bound. Similarly,
\begin{align*}
0 & =\max\left\{ GX,IGH^{\prime},CI^{\prime},X^{\prime}CH\right\} ,\\
 & =\max\left\{ I,I^{\prime}\right\} ,
\end{align*}
does not have a solution and so confirms $GC$ as a lower bound. We
still need to prove that $G$ or $C$ are not lower bounds in both
cases. So let $G=1$ then 
\begin{align*}
0 & =\max\left\{ HIX^{\prime},GX,IGH^{\prime},CI^{\prime}\right\} ,\\
 & =\max\left\{ HIX^{\prime},IH^{\prime},CI^{\prime}\right\} ,
\end{align*}
has $C=I=0$ as solution. For $C=1$ we get 
\begin{align*}
0 & =\max\left\{ HIX^{\prime},GX,IGH^{\prime},CI^{\prime}\right\} ,\\
 & =\max\left\{ HX^{\prime},GX,GH^{\prime}\right\} ,
\end{align*}
has $H=G=0$ as solution. In the remaining case 
\[
0=\max\left\{ GX,IGH^{\prime},CI^{\prime},X^{\prime}CH\right\} ,
\]
we see that for $G=1$ we get 
\[
0=\max\left\{ IH^{\prime},CI^{\prime},CH\right\} 
\]
has $C=H^{\prime}=0$ as solution. For $C=1$ we get 
\[
0=\max\left\{ GX,GH^{\prime},X^{\prime}H\right\} 
\]
has $H=G=0$ as solution. So $GC$ is really the largest lower bound
in both cases. That there are actually two solutions makes the problem
wrongly stated, as J. C.Wilson found out after Lewis Carroll sent
the challenge to him.

Since we found $IG=0$, we could equivalently simplify (3) to get
the logically equivalent assumption 
\begin{multicols}{2}

\begin{enumerate}
\item $HIX^{\prime}=0$ 
\item $GX=0$, 
\item $IG=0$, 
\item $CI^{\prime}=0$, 
\item $X^{\prime}CH=0$.\\
 
\end{enumerate}
\end{multicols}
 Here, however, we have $G,C,H$ as retinends and $X,I$ as eliminands!
Let $GCH=1$ then 
\begin{align*}
0 & =\max\left\{ HIX^{\prime},GX,IG,CI^{\prime},X^{\prime}CH\right\} ,\\
 & =\max\left\{ X,I,I^{\prime},X^{\prime}\right\} ,
\end{align*}
which is obviously false and so confirms that $GCH$ is a lower bound.
Are $GC,GH,CH$ still lower bounds? 
\begin{casenv}
\item $GC=1$ 
\begin{align*}
0 & =\max\left\{ HIX^{\prime},GX,IG,CI^{\prime},X^{\prime}CH\right\} ,\\
 & =\max\left\{ HIX^{\prime},X,I,I^{\prime},X^{\prime}H\right\} ,
\end{align*}
\item $GH=1$ 
\begin{align*}
0 & =\max\left\{ IX^{\prime},X,I,CI^{\prime},X^{\prime}C\right\} ,\\
 & =\max\left\{ X,I,C\right\} ,
\end{align*}
\item $CH=1$ 
\begin{align*}
0 & =\max\left\{ HIX^{\prime},GX,IG,CI^{\prime},X^{\prime}CH\right\} ,\\
 & =\max\left\{ G\right\} .
\end{align*}
\end{casenv}
We see that $GC$ is still the only largest lower bound and so 
\[
GC=0
\]
is still the proper complete conclusion. If we wanted to remember
the non-emptiness constraint $C\neq0$ at this point, we would have
to rephrase that as `all cyclists fail to be grocers'.

\newpage
\subsection{Some More Applications}

\subsubsection{Fat Lawyers}

The problem in Book XIV, Problems and Exercises 49, \cite{book834284},
has the following list of assumptions:

(1) No fat old men are active,

(2) All wise men are old,

(3) No unwise lawyers are fat.

Taking the universe $\mathcal{U}="men"$ and $a="active"$, $b="fat"$,
$c="lawyers"$, $d="old"$, $e="wise"$, we obtain the initial translation
\begin{multicols}{3}

\begin{enumerate}
\item $bd\subseteq a^{\prime}$, 
\item $e\subseteq d$, 
\item $e^{\prime}c\subseteq b^{\prime}$. 
\end{enumerate}
\end{multicols}
Lewis Carroll's version as a system of equations (nullities): 
\begin{multicols}{3}

\begin{enumerate}
\item $bda=0,$ 
\item $ed^{\prime}=0$, 
\item $e^{\prime}cb=0$. 
\end{enumerate}
\end{multicols}

We see that $e,d$ are eliminands and so $a,b,c$ are retinends.

So, we start with 
\[
abc
\]
as a possible lower bound of $\max\left\{ abd,ed^{\prime},bce^{\prime}\right\} $.
So let $abc=1$ then 
\begin{align*}
0 & =\max\left\{ abd,ed^{\prime},bce^{\prime}\right\} ,\\
 & =\max\left\{ d,ed^{\prime},e^{\prime}\right\} ,
\end{align*}
does not have a solution, which confirms that $abc$ is a lower bound.
Larger lower bounds can be excluded along the lines of earlier reasoning.
Thus we obtain that $abc=0$ is the proper complete conclusion, i.e.
``no fat lawyer is active''.

\subsubsection{Creatures in the Garden }

This problem is in Book XIV, Problems and Exercises 51, \cite{book834284}.

With $\mathcal{U}="creatures"$, $a="dignified"$, $b="greedily\:eating"$,
$c="fat"$,$d="frogs"$, $e="in\:this\:garden"$ we transcribe

(1) No frogs are dignified,

(2) All fat creatures in this garden eat greedily,

(3) Nothing that eats greedily is fat, except a frog,

first to 
\begin{multicols}{3}

\begin{enumerate}
\item $d\subseteq a^{\prime}$, 
\item $ce\subseteq b$, 
\item $b^{\prime}d^{\prime}\subseteq c$ 
\end{enumerate}
\end{multicols}

and then to nullities (again we ignore statements of non-emptiness)
\begin{multicols}{3}

\begin{enumerate}
\item $da=0$, 
\item $ceb^{\prime}=0$, 
\item $b^{\prime}d^{\prime}c^{\prime}=0$. 
\end{enumerate}
\end{multicols}

$c,d$ are eliminands, $a,b^{\prime},e$ are retinends. We suspect
\[
\max\left\{ ad,b^{\prime}ce,b^{\prime}c^{\prime}d^{\prime}\right\} \geq ab^{\prime}e.
\]
Indeed, if $ab^{\prime}e=1$ we have that 
\begin{align*}
0 & =\max\left\{ ad,b^{\prime}ce,b^{\prime}c^{\prime}d^{\prime}\right\} ,\\
 & =\max\left\{ d,c,c^{\prime}d^{\prime}\right\} ,
\end{align*}
has no solution. We check if $ab^{\prime}$, $ae$, $b^{\prime}e$
are still lower bounds. So, let $ab^{\prime}=1$ 
\begin{align*}
0 & =\max\left\{ ad,b^{\prime}ce,b^{\prime}c^{\prime}d^{\prime}\right\} ,\\
 & =e.
\end{align*}
For $ae=1$ 
\begin{align*}
0 & =\max\left\{ ad,b^{\prime}ce,b^{\prime}c^{\prime}d^{\prime}\right\} ,\\
 & =b^{\prime}.
\end{align*}
For $b^{\prime}e=1$ 
\begin{align*}
0 & =\max\left\{ ad,b^{\prime}ce,b^{\prime}c^{\prime}d^{\prime}\right\} ,\\
 & =a.
\end{align*}
In all three case we see that they are not lower bounds anymore. Thus
$ab^{\prime}e$ is the largest lower bound and so 
\[
ab^{\prime}e=0
\]
the proper complete conclusion, which reads `dignified creatures in
my garden are eating greedily'.

The next example we want to use as an illustration involves finding
the right universe $\mathcal{U}$, which is a bit more complex than
usually.

\subsubsection{The Problem of the Schoolboys}

\textquotedbl All the boys, in a certain school, sit together in one
large room every evening. They are of no less than five nationalities
- English, Scotch, Welsh, Irish, and German. One of the Monitors (who
is a great reader of Wilkie Collins' novels) is very observant and
takes MS. notes of almost everything that happens, with the view of
being a good sensational witness, in case any conspiracy to commit
a murder should be afoot. The following are some of his notes:

(1) Whenever some of the English boys are singing 'Rule, Britannia,'
and some not, some of the Monitors are wide awake;

(2) Whenever some of the Scotch are dancing reels, and some of the
Irish fighting, some of the Welsh are eating toasted cheese;

(3) Whenever all the Germans are playing chess, some of the Eleven
are not oiling their bats;

(4) Whenever some of the Monitors are asleep, and some not, some of
the Irish are fighting.

(5) Whenever some of the Germans are playing chess, and none of the
Scotch are dancing reels, some of the Welsh are not eating toasted
cheese.

(6) Whenever some of the Scotch are not dancing reels, and some of
the Irish are not fighting, some of the Germans are playing chess;

(7) Whenever some of the Monitors are awake, and some of the Welsh
are eating toasted cheese, none of the Scotch are dancing reels;

(8) Whenever some of the Germans are not playing chess, and some of
the Welsh are not eating toasted cheese, none of the Irish are fighting;

(9) Whenever all of the English are singing 'Rule, Britannia,\textquotedbl{}
and some of the Scotch are not dancing reels, none of the Germans
are playing chess;

(10) Whenever some of the English are singing \textquotedbl Rule,
Britannia,\textquotedbl , and some of the Monitors are asleep, some
of the Irish are not fighting;

(11) Whenever some of the Monitors are awake, and some of the Eleven
are not oiling their bats, some of the Scotch are dancing reels;

(12) Whenever some of the English are singing \textquotedbl Rule,
Britannia,' and some of the Scotch are not dancing reels, ...

Here the MS. breaks off suddenly. The problem is to complete the sentence,
if possible.\textquotedbl{}

With $\mathcal{V}=$''people at a boarding school'', we have a particular
interest in certain subsets labelled $X=\left(a,b,c,d,e,f,g,h,j,k,l,m,n,p\right)$:

$a=$''English boys'', $b=$''those singing 'Rule, Britannia'~``,
$c=$''Monitors'', $d=$''those asleep'', $e=$''Scottish boys'',
$f=$''those dancing reels'', $g=$''Irish boys'', $h=$''those
fighting'', $j=$''Welsh boys'', $k=$''those eating toasted cheese'',
$l=$''members of the Eleven'', $m=$''those oiling their bats'',
$n=$''German boys'', $p=$''those playing chess''. The actual
premises are now propositions about predicates $A_{k}\in2^{\left(2^{\mathcal{V}}\right)^{14}}$,
$k=0,\ldots,13$, about pairs of these subsets. This makes our universe
$\mathcal{U}=2^{\left(2^{\mathcal{V}}\right)^{14}}$. Our predicates
of interest are 
\begin{align*}
A_{0}\left(X\right) & =\left(0\neq ab\right), & A_{5}\left(X\right) & =\left(0\neq ef^{\prime}\right), & A_{10}\left(X\right) & =\left(0\neq lm\right),\\
A_{1}\left(X\right) & =\left(0\neq ab^{\prime}\right), & A_{6}\left(X\right) & =\left(0\neq gh\right), & A_{11}\left(X\right) & =\left(0\neq lm^{\prime}\right),\\
A_{2}\left(X\right) & =\left(0\neq cd\right), & A_{7}\left(X\right) & =\left(0\neq gh^{\prime}\right), & A_{12}\left(X\right) & =\left(0\neq np\right),\\
A_{3}\left(X\right) & =\left(0\neq cd^{\prime}\right), & A_{8}\left(X\right) & =\left(0\neq jk\right), & A_{13}\left(X\right) & =\left(0\neq np^{\prime}\right).\\
A_{4}\left(X\right) & =\left(0\neq ef\right), & A_{9}\left(X\right) & =\left(0\neq jk^{\prime}\right),
\end{align*}
Here $A_{10}$ is only listed for systematic reasons and is not actually
used in the premises.

The $A_{k},k=0,\ldots,13,$ are in $\mathcal{U}=2^{\left(2^{\mathcal{V}}\right)^{14}}$and
so the premises are in $2^{\mathcal{U}}$and can now be formalised
in Lewis Carroll's normal form as 
\begin{multicols}{2}

\begin{enumerate}
\item $A_{0}A_{1}A_{3}^{\prime}=0$, 
\item $A_{4}A_{6}A_{8}^{\prime}=0$, 
\item $A_{13}^{\prime}A_{11}^{\prime}=0$, 
\item $A_{2}A_{3}A_{6}^{\prime}=0$, 
\item $A_{12}A_{4}^{\prime}A_{9}^{\prime}=0$, 
\item $A_{5}A_{7}A_{12}^{\prime}=0$, 
\item $A_{3}A_{8}A_{4}=0$, 
\item $A_{13}A_{9}A_{6}=0$, 
\item $A_{1}^{\prime}A_{5}A_{12}=0$, 
\item $A_{0}A_{2}A_{7}^{\prime}=0$, 
\item $A_{3}A_{11}A_{4}^{\prime}=0$.\\
 
\end{enumerate}
\end{multicols}
 with retinends $A_{0},A_{2},A_{5}$ and eliminands $A_{1},A_{3},A_{4},A_{6},A_{7},A_{8},A_{9},A_{11},A_{12}$, $A_{13}$.
Thus, 
\[
A_{0}A_{2}A_{5}
\]
is our candidate for a lower bound. So let us assume that $A_{0}A_{2}A_{5}=1$
and then let us see if we can have 
\begin{align*}
0 & =\max\left\{ A_{0}A_{1}A_{3}^{\prime},A_{4}A_{6}A_{8}^{\prime},A_{13}^{\prime}A_{11}^{\prime},A_{2}A_{3}A_{6}^{\prime},A_{12}A_{4}^{\prime}A_{9}^{\prime},A_{5}A_{7}A_{12}^{\prime},\right.\\
 & \hspace*{1.3cm}\left.A_{3}A_{8}A_{4},A_{13}A_{9}A_{6},A_{1}^{\prime}A_{5}A_{12},A_{0}A_{2}A_{7}^{\prime},A_{3}A_{11}A_{4}^{\prime}\right\} \\
 & =\max\left\{ A_{1}A_{3}^{\prime},A_{4}A_{6}A_{8}^{\prime},A_{13}^{\prime}A_{11}^{\prime},A_{3}A_{6}^{\prime},A_{12}A_{4}^{\prime}A_{9}^{\prime},A_{7}A_{12}^{\prime},A_{3}A_{8}A_{4},\right.\\
 & \hspace*{1.3cm}\left.A_{13}A_{9}A_{6},A_{1}^{\prime}A_{12},A_{7}^{\prime},A_{3}A_{11}A_{4}^{\prime}\right\} \\
 & =\max\left\{ A_{1}A_{3}^{\prime},A_{4}A_{6}A_{8}^{\prime},A_{13}^{\prime}A_{11}^{\prime},A_{3}A_{6}^{\prime},A_{12}A_{4}^{\prime}A_{9}^{\prime},A_{12}^{\prime},A_{3}A_{8}A_{4},\right.\\
 & \hspace*{1.3cm}\left.A_{13}A_{9}A_{6},A_{1}^{\prime}A_{12},A_{3}A_{11}A_{4}^{\prime}\right\} \\
 & =\max\left\{ A_{3}^{\prime},A_{4}A_{6}A_{8}^{\prime},A_{13}^{\prime}A_{11}^{\prime},A_{3}A_{6}^{\prime},A_{4}^{\prime}A_{9}^{\prime},A_{3}A_{8}A_{4},A_{13}A_{9}A_{6},A_{3}A_{11}A_{4}^{\prime}\right\} \\
 & =\max\left\{ A_{4}A_{6}A_{8}^{\prime},A_{13}^{\prime}A_{11}^{\prime},A_{6}^{\prime},A_{4}^{\prime}A_{9}^{\prime},A_{8}A_{4},A_{13}A_{9}A_{6},A_{11}A_{4}^{\prime}\right\} \\
 & =\max\left\{ A_{4}A_{8}^{\prime},A_{13}^{\prime}A_{11}^{\prime},A_{4}^{\prime}A_{9}^{\prime},A_{8}A_{4},A_{13}A_{9},A_{11}A_{4}^{\prime}\right\} \\
 & =\begin{cases}
\max\left\{ A_{13}^{\prime}A_{11}^{\prime},A_{9}^{\prime},A_{13}A_{9},A_{11}\right\}  & \text{for }A_{4}=0,\\
\max\left\{ A_{8}^{\prime},A_{13}^{\prime}A_{11}^{\prime},A_{8},A_{13}A_{9}\right\}  & \text{for }A_{4}=1
\end{cases}\\
 & =\begin{cases}
\max\left\{ A_{13}^{\prime},A_{13}\right\}  & \text{for }A_{4}=0,\\
\max\left\{ A_{8}^{\prime},A_{8}\right\}  & \text{for }A_{4}=1
\end{cases}
\end{align*}

which clearly has no solution, confirming that 
\[
A_{0}A_{2}A_{5}=0
\]
or equivalently 
\[
A_{0}A_{5}\subseteq A_{2}^{\prime}
\]
is implied by the assumptions. Thus, the missing part in the final
statement (12) is 
\[
A_{2}^{\prime}\left(X\right)\left(\Leftrightarrow cd=0\Leftrightarrow c\subseteq d^{\prime}\right),
\]
which can be verbalised as ``none of the Monitors are asleep''.
Thus, the final statement (12) says: ``Whenever some of the English
are singing ``Rule, Britannia,'' and some of the Scotch are not
dancing reels, none of the Monitors are asleep''.

For more examples see Lewis Carroll, Symbolic Logic, \cite{Carroll1958,book834284}.
We conclude with a problem posed by Lewis Carroll, referred to as
``The Froggy Problem'', \cite[Book XIII, Chaper V]{book834284}),
which leads to an interesting new perspective on possible universes.

\subsubsection{The Froggy Problem}

The verbal premises are:

(1) When the day is fine, I tell Froggy ``You're quite the dandy,
old chap!'';

(2) Whenever I let Froggy forget that $\pounds$10 he owes me, and he begins
to strut about like a peacock, his mother declares ``He shall not
go out a-wooing!'';

(3) \textsl{Now} that Froggy's hair is out of curl, he has put away
his gorgeous waistcoat;

(4) Whenever I go out on the roof to enjoy a quiet cigar, I'm sure
to discover that my purse is empty;

(5) When my tailor calls with his little bill, and I remind Froggy
of that $\pounds$10 he owes me, he does not grin like a hyena;

(6) When it is very hot, the thermometer is high;

(7) When the day is fine, and I'm not in the humor for a cigar, and
Froggy is grinning like a hyena, I never venture to hint that he's
quite the dandy;

(8) When my tailor calls with his little bill and finds me with an
empty purse, I remind Froggy of that $\pounds$10 he owes me;

(9) \textsl{(Now)} My railway-shares are going up like anything!

(10) When my purse is empty, and when, noticing that Froggy has got
his gorgeous waistcoat on, I venture to remind him of that $\pounds$10 he
owes me, things are apt to get rather warm;

(11) \textsl{Now} that it looks like rain, and Froggy is grinning
like a hyena, I can do without my cigar;

(12) When the thermometer is high, you need not trouble yourself to
take an umbrella;

(13) When Froggy has his gorgeous waistcoat on, but is not strutting
about like a peacock, I betake myself to a quiet cigar;

(14) When I tell Froggy that he's quite the dandy, he grins like
a hyena;

(15) When my purse is tolerably full, and Froggy's hair is one mass
of curls, and when he is not strutting about like a peacock, I go
out on the roof;

(16) When my railway-shares are going up, and when it is chilly and
looks like rain, I have a quiet cigar;

(17) When Froggy's mother lets him go a-wooing, he seems nearly mad
with joy, and puts on a waistcoat that is gorgeous beyond words;

(18) When it is going to rain, and I am having a quiet cigar, and
Froggy is not intending to go a-wooing, you had better take an umbrella;

(19) When my railway-shares are going up, and Froggy seems nearly
mad with joy, that is the time my tailor always chooses for calling
with his little bill;

(20) When the day is cool and the thermometer low, and I say nothing
to Froggy about his being quite the dandy, and there's not the ghost
of a grin on his face, I haven't the heart for my cigar!

We have emphasised the `now' in the above to draw attention to the
fact that time has now to be included in our universe (no pun intended).
Lewis Carroll speaks about cosmophases to incorporate temporal statements.
So, `now' is phase of the cosmos as much as the others for example
the phase when `Froggy's hair is out of curl'. Recall, however, that
we can always construct a suitable universe be taking the union of
all `phases' occurring in the setup without having a nice common label.

Let us choose to follow Lewis Carroll's suggestion to take $\mathcal{U}$=''Cosmophases'',
and use the labels

$\epsilon$ = ``now'';

a = ``Froggy's hair is out of curl'';

b = ``Froggy intends to go a-wooing'';

c = ``Froggy is grinning like a hyena'';

d = ``Froggy's mother permits him to go a-wooing'';

e = ``Froggy seems nearly mad with joy'';

h = ``Froggy is strutting about like a peacock'';

l = ``I go out on my roof'';

k = ``Froggy is wearing a waistcoat that is gorgeous beyond words'';

m = ``I remind Froggy of the $\pounds$10 he owes me'';

n = ``I take a quiet cigar'';

r = ``I tell Froggy that he's quite the dandy'';

s = ``It is going to rain'',

t = ``It is very hot'';

v = ``My purse is empty'';

w = ``My railway shares are going up'';

z = ``My tailor calls with his little bill'';

A = ``The thermometer is high'';

B = ``You had better take an umbrella'';

The name for the cosmophase $\epsilon$, however, is not a proposition. A corresponding proposition
would be ``It is now''. The cosmophases $s$ and $B$ are also not
described as a corresponding proposition. Cosmophase s = ``It is
going to rain'' should correspond to a proposition such as ``It
looks like it is going to rain'' (otherwise we would assign truth
values to something that is going to happen, which is a prediction
rather than a proposition) and cosmophase B = ``You had better take
an umbrella'' should correspond to a proposition such as ``I recommend
you'd better take an umbrella''. These liberties in wording are typical
of Lewis Carroll and should not lead to much confusion. Indeed, let
us not be pedantic about this, we ``know'' what was meant.

A first translation into formal language yields 
\begin{multicols}{3}

\begin{enumerate}
\item $s^{\prime}\subseteq r$, 
\item $m^{\prime}h\subseteq d^{\prime}$, 
\item $\epsilon\subseteq ak^{\prime}$, 
\item $ln\subseteq v$, 
\item $zm\subseteq c^{\prime}$, 
\item $t\subseteq A$, 
\item $s^{\prime}n^{\prime}c\subseteq r^{\prime}$, 
\item $zv\subseteq m$, 
\item $\epsilon\subseteq w$, 
\item $vkm\subseteq t$, 
\item $\epsilon\subseteq scn^{\prime}$, 
\item $A\subseteq B^{\prime}$, 
\item $kh^{\prime}\subseteq n$, 
\item $r\subseteq c$, 
\item $v^{\prime}a^{\prime}h^{\prime}\subseteq l$, 
\item $wt^{\prime}s\subseteq n$, 
\item $d\subseteq ek$, 
\item $snb^{\prime}\subseteq B$, 
\item $we\subseteq z$, 
\item $t^{\prime}A^{\prime}r^{\prime}c^{\prime}\subseteq n^{\prime}$.\\
 
\end{enumerate}
\end{multicols}
 In terms of nullities we have 19 variables and 24 equations 
\begin{multicols}{3}

\begin{enumerate}
\item $s^{\prime}r^{\prime}=0$, 
\item $m^{\prime}hd=0$, 
\item[(3a)] $\epsilon a^{\prime}=0$, 
\item[(3b)] $\epsilon k=0$,\addtocounter{enumi}{1} 
\item $lnv^{\prime}=0$, 
\item $zmc=0$, 
\item $tA^{\prime}=0$, 
\item $s^{\prime}n^{\prime}cr=0$, 
\item $zvm^{\prime}=0$, 
\item $\epsilon w^{\prime}=0$, 
\item $vkmt^{\prime}=0$, 
\item[(11a)] $\epsilon s^{\prime}=0$, 
\item[(11b)] $\epsilon c^{\prime}=0$, 
\item[(11c)] $\epsilon n=0$,\addtocounter{enumi}{1} 
\item $AB=0$, 
\item $kh^{\prime}n^{\prime}=0$, 
\item $rc^{\prime}=0$, 
\item $v^{\prime}a^{\prime}h^{\prime}l^{\prime}=0$, 
\item $wt^{\prime}sn^{\prime}=0$, 
\item[(17a)] $de^{\prime}=0$, 
\item[(17b)] $dk^{\prime}=0$,\addtocounter{enumi}{1} 
\item $snb^{\prime}B^{\prime}=0$, 
\item $wez^{\prime}=0$, 
\item $t^{\prime}A^{\prime}r^{\prime}c^{\prime}n=0$. 
\end{enumerate}
\end{multicols}

Here the eliminands are $c,e,h,k,l,m,n,r,s,t,v,w,z,A,B$ and correspondingly
the retinends are $\epsilon,a^{\prime},b^{\prime},d.$

To check if $\epsilon a^{\prime}b^{\prime}d$ is a lower bound for
the maximum of the 24 products in the above system of nullities, we
assume $\epsilon a^{\prime}b^{\prime}d=1$ and see how that impacts
\begin{align*}
\mathcal{A}\coloneqq\max\left\{ s^{\prime}r^{\prime},m^{\prime}hd,\epsilon a^{\prime},\epsilon k,lnv^{\prime},zmc,tA^{\prime},s^{\prime}n^{\prime}cr,zvm^{\prime},\epsilon w^{\prime},vkmt^{\prime},\epsilon s^{\prime},\epsilon c^{\prime},\right.\\
\left.\epsilon n,AB,kh^{\prime}n^{\prime},rc^{\prime},v^{\prime}a^{\prime}h^{\prime}l^{\prime},wt^{\prime}sn^{\prime},de^{\prime},dk^{\prime},snb^{\prime}B^{\prime},wez^{\prime},t^{\prime}A^{\prime}r^{\prime}c^{\prime}n\right\} .
\end{align*}
Since $\epsilon a^{\prime}$ is part of the assumption, we already
have that $\epsilon a^{\prime}$ is a largest lower bound. In comparison,
the initial $\epsilon a^{\prime}b^{\prime}d$ is just a smaller\footnote{Via Lewis Carroll's own method Graham Hawker has derived that $\epsilon a^{\prime}b^{\prime}d=0$
can be obtained from the assumptions, \cite[Problems and Exercises - 30]{GH1}.
This approach is, however, much more involved than the method proposed
here. Moreover, here we are looking for not just any conclusion but for the sharpest one. All weaker conclusions are essentially obtainable as a combinatorial exercise from the sharpest one by taking products with other factors appearing in the collection of retinends.} non-trivial lower bound.Since $\epsilon a^{\prime}$ 
is a part of the assumption, we would want to discard this, leaving
us to test if there are other lower bounds like $\epsilon b'd$ or
$a'b'd$. If $\epsilon b'd=1$ we get to see if 
\begin{align*}
0 & =\max\left\{ s^{\prime}r^{\prime},m^{\prime}h,a^{\prime},k,lnv^{\prime},zmc,tA^{\prime},s^{\prime}n^{\prime}cr,zvm^{\prime},w^{\prime},vkmt^{\prime},s^{\prime},c^{\prime},\right.\\
 & \hspace*{1.3cm}\left.n,AB,kh^{\prime}n^{\prime},rc^{\prime},v^{\prime}a^{\prime}h^{\prime}l^{\prime},wt^{\prime}sn^{\prime},e^{\prime},k^{\prime},snB^{\prime},wez^{\prime},t^{\prime}A^{\prime}r^{\prime}c^{\prime}n\right\} \\
 & =\max\left\{ m^{\prime}h,a^{\prime},zm,tA^{\prime},zvm^{\prime},w^{\prime},AB,v^{\prime}a^{\prime}h^{\prime}l^{\prime},wt^{\prime},e^{\prime},1,wez^{\prime}\right\} ,
\end{align*}
has a solution. This is not the case, showing that $\epsilon b^{\prime}d$
is a lower bound. Collect all 2-factor products and assume $\epsilon d=1$
\begin{multline*}
\max\left\{ s^{\prime}r^{\prime},\epsilon a^{\prime},\epsilon k,tA^{\prime},\epsilon w^{\prime},\epsilon s^{\prime},\epsilon c^{\prime},\epsilon n,AB,rc^{\prime},de^{\prime},dk^{\prime}\right\} \\
=\max\left\{ s^{\prime}r^{\prime},a^{\prime},k,tA^{\prime},w^{\prime},s^{\prime},c^{\prime},n,AB,rc^{\prime},e^{\prime},k^{\prime}\right\} =1.
\end{multline*}
Thus, $\epsilon d$ is already a largest lower bound for this part,
found as a lucky guess, the more so for the whole\footnote{That this actually works points at massive redundancies in the assumption.
These redundancies will help us later to `guess' solutions to nullities
in order to exclude other candidates for lower bounds. } $\mathcal{A}$. In particular, the assumption $\mathcal{A}=0$ implies
\begin{equation}
\epsilon d=0,\label{eq:Froggy-solve}
\end{equation}
which may be verbalised as `Froggy's mother does not permit him to
go a-wooing now'. This fits in with the first verse of the poem \textquotedbl A
Frog He Would A-wooing Go\textquotedbl : ``A frog he would a-wooing
go, Whether his mother would let him or no.''\footnote{Here is a complete text of this old song:\\
 A Frog he would a-wooing go, Heigh-ho, says Rowley, A Frog he would
a-wooing go, Whether his mother would let him or no, With a Roley,
Poley, Gammon and Spinach, Heigh-ho says Anthony Rowley.

He saddled and bridled a great black snail, Heigh-ho, says Rowley,
He saddled and bridled a great black snail, And rode between the horns
and the tail, With a Roley, Poley, Gammon and Spinach, Heigh-ho says
Anthony Rowley.

So off he set with his opera hat, Heigh-ho, says Rowley, So off he
set with his opera hat, And on the way he met with a rat, With a Roley,
Poley, Gammon and Spinach, Heigh-ho says Anthony Rowley.

They rode till they came to Mousey Hall, Heigh-ho, says Rowley, They
rode till they came to Mousey Hall, And there they both did knock
and call, With a Roley, Poley, Gammon and Spinach, Heigh-ho says Anthony
Rowley.

\textquotedbl Pray, Mrs. Mouse, are you within?\textquotedbl{} Heigh-ho,
says Rowley, \textquotedbl Pray, Mrs. Mouse, are you within?\textquotedbl{}
\textquotedbl Oh yes, sir, here I sit and spin.\textquotedbl{} With
a Roley, Poley, Gammon and Spinach, Heigh-ho says Anthony Rowley.

Then Mrs. Mouse she did come down, Heigh-ho, says Rowley, Then Mrs.
Mouse she did come down, All smartly dressed in a russet gown, With
a Roley, Poley, Gammon and Spinach, Heigh-ho says Anthony Rowley.

\textquotedbl Pray, Mrs. Mouse, can you give us some beer,\textquotedbl{}
Heigh-ho, says Rowley, \textquotedbl Pray, Mrs. Mouse, can you give
us some beer, That froggy and I may have good cheer?\textquotedbl{}
With a Roley, Poley, Gammon and Spinach, Heigh-ho says Anthony Rowley.

She had not been sitting long to spin, Heigh-ho, says Rowley, She
had not been sitting long to spin, When the cat and the kittens came
tumbling in, With a Roley, Poley, Gammon and Spinach, Heigh-ho says
Anthony Rowley.

The cat she seized Master Rat by the crown, Heigh-ho, says Rowley,
The cat she seized Master Rat by the crown, The kitten she pulled
Miss Mousey down, With a Roley, Poley, Gammon and Spinach, Heigh-ho
says Anthony Rowley.

This put Mr. Frog in a terrible fright, Heigh-ho, says Rowley, This
put Mr. Frog in a terrible fright, He took up his hat and he wished
them \textquotedbl Good night!\textquotedbl{} With a Roley, Poley,
Gammon and Spinach, Heigh-ho says Anthony Rowley.

And as he was passing over the brook, Heigh-ho, says Rowley, And as
he was passing over the brook, A lily white duck came and gobbled
him up, With a Roley, Poley, Gammon and Spinach, Heigh-ho says Anthony
Rowley.

So there's an end of one, two, and three, Heigh-ho, says Rowley, So
there's an end of one, two, and three, The Rat, the Mouse, and little
Froggy, With a Roley, Poley, Gammon and Spinach, Heigh-ho says Anthony
Rowley.}

That $\epsilon d=0$ is actually the proper complete conclusion requires
the exclusions of other candidates. Excluding $\epsilon d$ and $\epsilon a^{\prime}$
from $\epsilon a^{\prime}b^{\prime}d$, we are left with $\epsilon b^{\prime}$
and $a^{\prime}b^{\prime}d$ as remaining candidates.

For $\epsilon b^{\prime}=1$ we look for solutions of 
\begin{align*}
0 & =\max\left\{ s^{\prime}r^{\prime},m^{\prime}hd,a^{\prime},k,lnv^{\prime},zmc,tA^{\prime},s^{\prime}n^{\prime}cr,zvm^{\prime},w^{\prime},vkmt^{\prime},s^{\prime},c^{\prime},n,\right.\\
 & \hspace*{1.3cm}\left.AB,kh^{\prime}n^{\prime},rc^{\prime},v^{\prime}a^{\prime}h^{\prime}l^{\prime},wt^{\prime}sn^{\prime},de^{\prime},dk^{\prime},snB^{\prime},wez^{\prime},t^{\prime}A^{\prime}r^{\prime}c^{\prime}n\right\} \\
 & =\max\left\{ zm,zvm^{\prime},B,ez^{\prime}\right\} ,
\end{align*}
which is for example given by 
\[
z=B=e=0.
\]
This shows that $\epsilon b^{\prime}$ is no lower bound. So, let
finally $a^{\prime}b^{\prime}d=1$ then 
\begin{align*}
0 & =\max\left\{ s^{\prime}r^{\prime},m^{\prime}hd,\epsilon,\epsilon k,lnv^{\prime},zmc,tA^{\prime},s^{\prime}n^{\prime}cr,zvm^{\prime},\epsilon w^{\prime},vkmt^{\prime},\epsilon s^{\prime},\epsilon c^{\prime},\epsilon n,\right.\\
 & \hspace*{1.3cm}\left.AB,kh^{\prime}n^{\prime},rc^{\prime},v^{\prime}h^{\prime}l^{\prime},wt^{\prime}sn^{\prime},e^{\prime},k^{\prime},snB^{\prime},wez^{\prime},t^{\prime}A^{\prime}r^{\prime}c^{\prime}n\right\} \\
 & =\max\left\{ s^{\prime}r^{\prime},m^{\prime}h,lnv^{\prime},zmc,tA^{\prime},s^{\prime}n^{\prime}cr,zvm^{\prime},vmt^{\prime},\right.\\
 & \hspace*{1.3cm}\left.AB,h^{\prime}n^{\prime},rc^{\prime},v^{\prime}h^{\prime}l^{\prime},wt^{\prime}sn^{\prime},snB^{\prime},wz^{\prime},t^{\prime}A^{\prime}r^{\prime}c^{\prime}n\right\} ,
\end{align*}
To chase for a zero, we try another guess by taking $w=z=n=h^{\prime}=m^{\prime}=0$
reducing our search to 
\[
0=\max\left\{ s^{\prime}r^{\prime},tA^{\prime},s^{\prime}n^{\prime}cr,vt^{\prime},AB,rc^{\prime}\right\} .
\]
There are still many literals in this equation, but apparently $s^{\prime}=B=A^{\prime}=t^{\prime}=c^{\prime}=0$
makes this indeed vanish. This confirms that $a^{\prime}b^{\prime}d$
is not a lower bound. With some more guesswork we can also include
that $\epsilon$ or $d$ separately are lower bounds, finally confirming
that \eqref{eq:Froggy-solve} is indeed the proper complete conclusion.
We try to substitute $s^{\prime}=B=A^{\prime}=t^{\prime}=c^{\prime}=w=z=n=h^{\prime}=m^{\prime}=0$
into $\mathcal{A}$ to obtain
\[
\max\left\{ \epsilon a^{\prime},\epsilon k,\epsilon,de^{\prime},dk^{\prime}\right\} .
\]
Now, if $d=1$ then we can still have
\begin{align*}
0 & =\max\left\{ \epsilon a^{\prime},\epsilon k,\epsilon,e^{\prime},k^{\prime}\right\} ,
\end{align*}
for $\epsilon=0$ and $e^{\prime}=k^{\prime}=0$. This ensures that
$d$ alone is not a lower bound for $\mathcal{A}$. For $\epsilon=1$
this would not lead to the desired outcome. Let us try $\epsilon=1$
and $d=s^{\prime}=B=A^{\prime}=t^{\prime}=c^{\prime}=w^{\prime}=n=h^{\prime}=m^{\prime}=0$\textasciiacute 
\begin{align*}
0 & =\max\left\{ a^{\prime},k,zmc,de^{\prime},dk^{\prime},ez^{\prime}\right\} ,\\
 & =\max\left\{ a^{\prime},k,zmc,ez^{\prime}\right\} .
\end{align*}
We read off that $a^{\prime}=k=z=e=0$ ensures that this equality
can hold, excluding also $\epsilon$ alone to become a lower bound.

\section{Computerised solutions}
Although this method can be carried out manually, it is even more effective with the help of a computer.
To do this, we have to encode the logical statements. As before, we will use nullities, 
which are essentially a set of logical variables that can be in one of three states: positive, negated, or unused.
We encode these states with the values \texttt{1}, \texttt{-1} and \texttt{0}, respectively. 
Let us look at an example of this process.

\vspace{1cm}

\textbf{The Pets Problem}
\begin{enumerate}
  \item All old monkeys are docile;
  \item All affectionate cats purr;
  \item All healthy dogs bark;
  \item All my fat pets are in the basket;
  \item No small pigs are affectionate;
  \item Black cats can be easily carried about;
  \item No cross active creatures jump, except monkeys;
  \item I make pets of all affectionate docile creatures;
  \item A healthy monkey, if electrified, becomes a pig;
  \item The large creatures in this basket are all cats;
  \item Old dogs are always cross;
  \item No fat creatures bark, except when electrified;
  \item All affectionate old creatures jump;
  \item All the creatures in this basket, that purr, are black;
  \item Large dogs cannot be easily carried about;
  \item Creatures that are cross and destitute of affection, are always lazy;
  \item Unhealthy creatures are either large or thin. 
\end{enumerate}

Find the complete conclusion, \cite[Chapter XIV]{book834284}.\\

Taking the universe 
$\mathcal{U}="\text{creatures}"$, 
and the variables

\begin{tabular}{lll}                 
  $a="\text{large}"$,                &
  $b="\text{in the basket}"$,        &
  $c="\text{cats}"$,                 \\
  $d="\text{docile}"$,               &
  $e="\text{easily carried about}"$, &
  $f="\text{affectionate}"$,         \\
  $g="\text{fat}"$,                  &
  $h="\text{healthy}"$,              &
  $j="\text{jump}"$,                 \\
  $k="\text{dogs}"$,                 &
  $l="\text{lazy}"$,                 &
  $m="\text{monkeys}"$,              \\
  $n="\text{black}"$,                &
  $o="\text{old}"$,                  &
  $p="\text{pets}"$,                 \\
  $r="\text{pigs}"$,                 &
  $s="\text{purring}"$,              &
  $t="\text{electrified}"$,          \\
  $w="\text{barking}"$,              &
  $x="\text{cross}"$,                
\end{tabular}

we obtain the translation

\begin{multicols}{2}
\begin{enumerate}
  \item $om\subseteq d \quad\Leftrightarrow\quad omd'=0$;
  \item $fc\subseteq s \quad\Leftrightarrow\quad fcs'=0$;
  \item $hk\subseteq w \quad\Leftrightarrow\quad hkw'=0$;
  \item $gp\subseteq b \quad\Leftrightarrow\quad gpb'=0$;
  \item $a^{\prime}r\subseteq f^{\prime} \quad\Leftrightarrow\quad a'rf=0$;
  \item $nc\subseteq e \quad\Leftrightarrow\quad nce'=0$;
  
  \item $xl^{\prime}m^{\prime}\subseteq j' \quad\Leftrightarrow\quad xl'm'j=0$;

  \item $fd\subseteq p \quad\Leftrightarrow\quad fdp'=0$;
  \item $hmt\subseteq r \quad\Leftrightarrow\quad hmtr'=0$;
  \item $ab\subseteq c \quad\Leftrightarrow\quad abc'=0$;
  \item $ok\subseteq x \quad\Leftrightarrow\quad okx'=0$;
  
  \item $gt^{\prime}\subseteq w^{\prime} \quad\Leftrightarrow\quad gt'w=0$;

  \item $of\subseteq j \quad\Leftrightarrow\quad ofj'=0$;
  \item $sb\subseteq n \quad\Leftrightarrow\quad sbn'=0$;
  \item $ak\subseteq e^{\prime} \quad\Leftrightarrow\quad ake=0$;
  \item $xf^{\prime}\subseteq l \quad\Leftrightarrow\quad xf'l'=0$;
  \item $h^{\prime}\subseteq\left(g^{\prime}\cup a\right) \Leftrightarrow\quad h'a'g=0$. 
\end{enumerate}
\end{multicols}

Thus, there are 
17 
nullities in 20 variables. Note that two statements are represented by two nullities each.
Each nullity now becomes one row of a 
$17\times 20$ 
matrix, with each column representing one variable. 

Note that we did not use the exclusive "or" statement in the last statement 17). 
It is unclear whether, in Lewis's work, ``'either \dots or'' would be exclusive without 
the additional 'but not both'. So, let us  first stick with the situation as stated.

\newpage
\begin{minipage}{\textwidth}
\begin{lstlisting}[%caption = {LewisCarrollExample.m},
%                  linerange={1-2,317-341,352-354}]{LewisCarrollExample.m}
                 linerange={1-2,337-359,524,531-532}]
function [A,names] = LewisCarrollExample(example)
  switch example
            0  1  0  0  0  1  0  0  0  0  0  0  0 -1  0  0;...
            1  0  0 -1  0  0  0  0  0  0  0  0  0  0  1  0;...
            0  0  0  0  0  0  1  0  0  0  1  0 -1  0  0  0;...
            0  0  0  1  0  0  0  0  0  0  0 -1  0  0  0  1;...
            0  0  0 -1  0  0  0  0  0  0  1  0  1  0  0  0;...
            1  0  0  0  0  0  0  0  0 -1  0  0  0  0  0 -1;...
            0  0  1 -1  0  0  1  0  0  0  0  0  0  0  0  0;...
            1  0  0  0  0  0  0  0 -1  0  0  1  0  0  0  0;...
            1  0  0  0  0  0  0  0  0  0 -1  1  0  0  0  0;...
            0  0  0  0  0  0  0  0  0 -1  1  0  0  0  1  0;...
            0  0  1  0 -1  0  0  0  1  0  0  0  0  0  0  1;...
            0 -1  0  0  0  0  1  0  0  0  0  0  0 -1  1  0;...
            1  0 -1  0  0  0  0  0  0  1  0  0  0  0  0  0;...
            0  0  0  1  0  0  0  1  0  0 -1  0  1  0  0  0;...
            0  0  0  0  1  0  0  0  0  0  1  1  0  0  0  0;...
            1  0  0  0  0  0  0  0  0  0  0  0  1  0  0 -1;...
            0  0  0  0  0 -1  0 -1  0  0  0  0  0  0  1  0;...
            0  0  0  1  0  0  0  0  0  1 -1 -1  0  0  0  0;...
          ];
      %     a  b  c  d  e  h  k  l  m  n  r  s  t  v  w  z
    case 31
      titletext = "Pets inclusive";
      names = {'a','b','c','d','e','f','g','h','j','k','l','m','n','o','p','r','s','t','w','x'};
      % a  b  c  d  e  f  g  h  j  k  l  m  n  o  p  r  s  t  w  x
      A = [...
        0  0  0  0  0  0  0  0  0  1  0  0  0  1  0  0  0  0  0 -1;...
        0  0  0  0  0  0  1  0  0  0  0  0  0  0  0  0  0 -1  1  0;...
        0  0  0  0  0  1  0  0 -1  0  0  0  0  1  0  0  0  0  0  0;...
function [A,names] = LewisCarrollExample(example)
  switch example
    case 31
      titletext = "Pets inclusive";
      names = {'a','b','c','d','e','f','g','h','j','k','l','m','n','o','p','r','s','t','w','x'};
      % a  b  c  d  e  f  g  h  j  k  l  m  n  o  p  r  s  t  w  x
      A = [...
        0  0  0 -1  0  0  0  0  0  0  0  1  0  1  0  0  0  0  0  0;...
        0  0  1  0  0  1  0  0  0  0  0  0  0  0  0  0 -1  0  0  0;...
        0  0  0  0  0  0  0  1  0  1  0  0  0  0  0  0  0  0 -1  0;...
        0 -1  0  0  0  0  1  0  0  0  0  0  0  0  1  0  0  0  0  0;...
       -1  0  0  0  0  1  0  0  0  0  0  0  0  0  0  1  0  0  0  0;... 
        0  0  1  0 -1  0  0  0  0  0  0  0  1  0  0  0  0  0  0  0;...
        0  0  0  0  0  0  0  0  1  0 -1 -1  0  0  0  0  0  0  0  1;...
        0  0  0  1  0  1  0  0  0  0  0  0  0  0 -1  0  0  0  0  0;...
        0  0  0  0  0  0  0  1  0  0  0  1  0  0  0 -1  0  1  0  0;...
        1  1 -1  0  0  0  0  0  0  0  0  0  0  0  0  0  0  0  0  0;...
        0  0  0  0  0  0  0  0  0  1  0  0  0  1  0  0  0  0  0 -1;...
        0  0  0  0  0  0  1  0  0  0  0  0  0  0  0  0  0 -1  1  0;...
        0  0  0  0  0  1  0  0 -1  0  0  0  0  1  0  0  0  0  0  0;...
        0  1  0  0  0  0  0  0  0  0  0  0 -1  0  0  0  1  0  0  0;...
        1  0  0  0  1  0  0  0  0  1  0  0  0  0  0  0  0  0  0  0;...
        0  0  0  0  0 -1  0  0  0  0 -1  0  0  0  0  0  0  0  0  1;...
       -1  0  0  0  0  0  1 -1  0  0  0  0  0  0  0  0  0  0  0  0;...
          ];
    otherwise
      % titletext = '';
      % names = {'','','','',''};
      % %     a  b  c  d  e
      % A = [...
      %       0  0  0  0  0;...
      %     ];
      error('example not implemented');
  end
  fprintf('Example %i: %s\n',example,titletext);
end
      error('example not implemented');
  end
\end{lstlisting}
\end{minipage} 

The result is the matrix \texttt{A}, which can be used in the actual solver. 
We will present it step by step. First, the data is loaded.

\begin{lstlisting}[%caption = {LewisCarrollLogic.m},
                 linerange={1-4}]
%% step 1) load example from list
example   = 31;
[A,names] = LewisCarrollExample(example);
[k,m]     = size(A); tic;
\end{lstlisting}

Now, \texttt{A} is a \texttt{k} by \texttt{m} matrix of the above type, and \texttt{names} contains the names of the variables. 
We identify the retinends and eliminands within the variables.

\begin{lstlisting}[%caption = {LewisCarrollLogic.m},
                 linerange={5-20},firstnumber=5]
%% step 2) find retinends and eliminands
id_pos     = any(A>0,1); id_neg = any(A<0,1);
eliminands = id_pos.*id_neg>0;
retinends  = ~eliminands;
r = sum(retinends);                    % number of retinends
A = [A(:,retinends) A(:,eliminands)];  % reorder A, retinends first
id_neg      = any(A(:,1:r)<0,1);       % which retinends are negated
A(:,id_neg) = -A(:,id_neg);            % positivity in A for retinends
fprintf('retinends:')                  % output
id = find(retinends);
for i=1:r
  if id_neg(i); fprintf(' not'); end
  fprintf(' (%s),',names{id(i)});
end;fprintf('\b\n');
fprintf('eliminants:');fprintf(' (%s),',names{eliminands});
fprintf('\b\n');
\end{lstlisting}

This is achieved by identifying columns containing entries of \texttt{-1} and \texttt{1}. 
The results are then presented as text.
Furthermore, \texttt{r} contains the number of retinends and \texttt{A} is reordered so that the first \texttt{r} columns correspond 
to the retinends and the remaining columns correspond to the eliminands. Furthermore, all entries for the retinends are now non-negative;
if an entry was previously negative, the negated variable is now used as a retinend. 

Now we come to the actual process of solving.

\begin{lstlisting}[%caption = {LewisCarrollLogic.m},
                 linerange={21-50},firstnumber=21]
%% step 3) find conclusions by eliminating retinends
conclusion = zeros(1,m); nr = 0;
for i=r-1:-1:0                         % #elminimated retinends
  subsets = nchoosek(1:r,i);           % subsets: i elements from 1:r
  for j=1:size(subsets,1)
    retelim = subsets(j,:);            % choose subset
    if r==1&&i==0; retelim = []; end   % treatment special case
    ret = setdiff(1:r,retelim);        % new set of retinends
    elim = [retelim,r+1:m];            %    and eleminants
    Ar = [A(:,ret) A(:,elim)];         % reorder A
    x = true(1,m); x(r+1-i:end)=false; % bitpattern of retinends
    y = bin2dec(num2str(x));           % binary representation
    sol = true;
    Arm = (Ar<=0); Arp = (Ar>=0);
    for z=0:2^(m-r)-1                  % bitpatterns start: 1s and 0s
      x = dec2bin(y+z)-'0';            % y+z has the pattern x
      C = any(all(((1-x).*Arm+...      % if Ar<0 we have ~x=1-x,
                       x.*Arp)>0,2));  % if Ar=0 we have 1=1-x+x, 
                                       % if Ar>0 we have x,
                                       % then all/min in rows
                                       % any/max  in cols
      sol = and(sol,C);                % if any x fails we can skip 
      if ~sol; break; end              %    these retinends
    end
    if sol                             % store conclusion found
      nr = nr+1;
      conclusion(nr,ret) = 1;
    end
  end
end
\end{lstlisting}

Solving is done iteratively for a given matrix $A=(a_{ij})\in\{-1,0,1\}^{k\times m}$ by the following algorithm:

\begin{quotation}
For each possible conclusion and corresponding selection of retinends, ask the following question: 
Does there exist an assignment $x\in 2^m$, such that all $x_r=1$ for all $r$ corresponding to the retinends, and such that
\[
  \max_i \min_j \left((a_{ij}\geq 0)x_j+(a_{ij}\leq 0)(1-x_j)\right)=0?
\]
If so, the considered conclusion is rejected.
\end{quotation}
Here the structure of the matrix $A$ is used. If $a_{ij}=1$, we use $x_j$ for the nullity, while for $a_{ij}=-1$ we use the negation $1-x_j$.
However, we also have to consider the case $a_{ij}=0$, in which the corresponding variable is not used. To apply the minimum
afterwards, we must take $1=x_j+(1-x_j)$ here, as this is not ``seen'' in the minimum. 

Note that this is a brute-force algorithmic solution and it considers all $2^r$ possible conclusions 
and $2^{m-i}$ possible assignments of $x$, where $i$ is the number of retinends chosen in this particular 
conclusion. This is feasible for small problem instances as presented so far.
The problem itself can actually be seen as a Boolean satisfiability problem (SAT), which is NP-hard and 
therefore generally only solvable at an exponential cost. 
However, there are solvers available that are quite fast in practise.

The result is a matrix \texttt{conclusion}, which stores the conclusions row-wise with \texttt{1} for a kept retinend in that conclusion.
These solutions may still contain superfluous conclusions, which will be removed in the next step.

\begin{minipage}[t]{\textwidth}  
\begin{lstlisting}[%caption = {LewisCarrollLogic.m},
                 linerange={51-64},firstnumber=51]
%% step 4) remove implications in the conclusion set
i = 0;
while i<nr
  i = i+1; j = i;
  id_i = find(conclusion(i,:));        % indices of conclusion_i
  while j<nr
    j = j+1;
    id_j = find(conclusion(j,:));      % indices of conclusion_j
    if ismember(id_i,id_j)             % id_i in id_j -> implication
      conclusion(j,:) = [];            % remove conclusion_j
      nr = nr-1; j = j-1;
    end
  end
end
\end{lstlisting}
\end{minipage} 

The first type of conclusion to be removed is that which is an implication of another conclusion.
The process of comparison is made easier by the storage method used for the possible conclusions. 
The ones with fewer variables are stored first, so it is only necessary to compare those with a 
smaller index with those with a larger index.
If one is a subset of another, the latter one is removed. Finally, to ensure proper conclusions, 
we also have to remove conclusions that are equivalent to assumptions.

\begin{lstlisting}[%caption = {LewisCarrollLogic.m},
                 linerange={65-85},firstnumber=65]
%% step 5) remove conclusions equivalent to assumptions
i = 0;
Am = (A<=0); Ap = (A>=0);
while i < nr
  i = i+1;
  set = false(k, 1);                   % logical mask
  for j = find(conclusion(i,:))        % iterate conclusion indices
    for z = 1:(2^m-1)                  % generate bit patterns 
      x = dec2bin(z,m)-'0';            %    on the fly
      if x(j); continue; end           % skip if x_j = 1
      C = all(((1-x).*Am+x.*Ap)>0, 2); % as above, but leave as vector
      set = set | C;                   % mark rows where C=Ax=1
      if all(set); break; end          % all rows => keep conclusion
    end
    if all(set); break; end            % no need to test other j
  end
  if ~all(set)                         % remove conclusion
    conclusion(i,:) = [];
    nr = nr-1; i = i-1;
  end
end
\end{lstlisting}

Once again, we opted for an inefficient implementation involving an exponential loop for simplicity.
As with the previous approach, this problem can be rewritten as a set of SATs, with one for each row 
of A and each possible conclusion.
Solving these problems is generally much faster.

Finally, we can output the proper conclusion.

\begin{minipage}[t]{\textwidth}  
\begin{lstlisting}[%caption = {LewisCarrollLogic.m},
                 linerange={86-95},firstnumber=86]
%% step 6) output conclusion set
fprintf('conclusion: ');
for i=1:nr
  ret = find(conclusion(i,:));
  fprintf('0 =');                
  for k=1:length(ret)
    if id_neg(ret(k)); fprintf(' not'); end
    fprintf(' (%s) and',names{id(ret(k))});
  end; fprintf('\b\b\b\n');
end; toc
\end{lstlisting}
\end{minipage} 

For the given example, the output of a run looks like this:

\begin{minipage}[t]{\textwidth}
\begin{lstlisting}
>> LewisCarrollLogic
Example 31: Pets inclusive
retinends: (g), (k), not (l), (o)
eliminands: (a), (b), (c), (d), (e), (f), (h), (j), (m), (n), (p), (r), (s), (t), (w), (x)
conclusion: 0 = (g) and (k) and not (l) and (o) 
\end{lstlisting}
\end{minipage} 

Thus, we get in natural language form ``old fat dogs are lazy''. 


Now let us add the nullity 
\[
  18)\quad h'ag'=0.
\]
to the assumptions effectively turning the inclusive "or" in 17) into an exclusive "or". 
Here the run yields:

\begin{minipage}[t]{\textwidth}
\begin{lstlisting}
>> LewisCarrollLogic
Example 33: Pets exclusive
retinends: (k), not (l), (o)
eliminants: (a), (b), (c), (d), (e), (f), (g), (h), (j), (m), (n), (p), (r), (s), (t), (w), (x)
conclusion: 
\end{lstlisting}
\end{minipage} 

Thus, this interpretation of the problem has only the trivial conclusion $0$. This suggests to dispose of the constraint 18).

\vspace{1cm}

\textbf{The London-Friends Problem}

An instance of the largest  logical problems we found in Lewis Carroll's work is in Problems and Exercises, Book XIV, Problem 1, \cite{book834284}, and goes like the following:
{\footnotesize
\begin{enumerate}
  \item I don't like walking with any London-friend of mine, unless he wears a tall hat;
  \item The Colonel is ready to play billiards with any man who is not the subject of general conversation, unless he happens to be fat;
  \item A man, who knows what o'clock it is and who never yields to me in argument, is unattractive;
  \item No old tight-rope-dancer ever rouses me to enthusiasm;
  \item Any novelist, whom I take pleasure in cutting, always cuts me dead;
  \item A man, who goes about in kid-gloves, but without his coat, is a humbug;
  \item All my intimate friends in London are young men;
  \item No man, who is the subject of general conversation, ever rouses me to enthusiasm, unless he is a horsey man;
  \item A man, who has his wits about him and does not choose his own wines, is always ``at home'' when I call;
  \item A man, who is a good shot and never tells pointless anecdotes, is sure to have a good temper;
  \item All humbugs, who write novels, are intimate friends of mine;
  \item I like to walk with a good-tempered man, unless he goes about in his shirt-sleeves;
  \item A man, who never loses his umbrella, and is not easily taken in, is sure to be an early riser;
  \item Fat men, who do not dance on the tight-rope, are universally respected;
  \item I regard with contemptuous pity a man who fails in life, and who runs across the street;
  \item A man, who does not stick to business, is not likely to be elected Mayor, unless he has bushy whiskers;
  \item An elephant hunter always rouses me to enthusiasm, unless he happens to be a farmer;
  \item Any London friend of mine, who tells pointless anecdotes, is a humbug;
  \item I never invite an old man to dinner, unless he has lent me money;
  \item A man, who does not stick to business, does not run across the street, and has bushy hair, is in no danger of getting a bad fall;
  \item A man, who gets up late and sometimes loses his umbrella, has little chance of marrying an heiress;
  \item An old man, who cares for appearances, always wears kid-gloves;
  \item A good-tempered man never cuts me dead, unless he is a humbug;
  \item A man, who never tells pointless anecdotes and has never lent me money, has his wits about him;
  \item A man, who chooses his own wines and always yields to me in argument, is the sort that I invite to dine with me;
  \item I always try to be civil to a man who fails in life, unless he has bushy whiskers;
  \item All farmers are horsey men;
  \item A novelist is a dull companion, unless he rouses me to enthusiasm;
  \item All men, who get up early, and stick to business and win universal respect, are rich;
  \item Any London friend of mine, to whom I try to be civil, will probably be elected Mayor;
  \item Any good-tempered man, who has lent me money and does not care for appearances, is willing to shake hands with me when I am in rags;
  \item The only men, with whom the Colonel will play billiards, are either horsey men or farmers;
  \item I always invite an attractive man to dine with me, provided he is rich;
  \item A man, who is apt to walk on tip-toe and whom I regard with contemptuous pity, is sure to be one who sticks to business;
  \item The only men, who are always ``at home'' to me, but whom I never invite to dinner, are magistrates;
  \item I always make an intimate friend of a man who will shake hands with me when I am in rags and will give up his umbrella to me when it is raining;
  \item Any London friend of mine, who understands horses, is universally respected;
  \item An unattractive man, who chooses his own wines, is easily taken in;
  \item I have sufficient courage to insult any novelist, unless he happens to be a good shot;
  \item An old man, who is apt to walk on tip-toe, will probably get a bad fall;
  \item A man, who never knows what o'clock it is, and who has never lent me money, will probably marry an heiress;
  \item No London-friend of mine, who has his wits about him, is easily taken in;
  \item I never forget any old man who is willing to shake hands with me when I am in rags;
  \item A novelist, who does not stick to business, is sure to fail in life;
  \item I do not dare to insult an ill-tempered man, unless he happens to be an intimate friend of mine;
  \item Those magistrates, who will not shake hands with me when I am in rags, always choose their own wine;
  \item All dull companions are either horsey men or elephant-hunters;
  \item Men, who wear tall hats and kid gloves, always cut me dead;
  \item A man, who has bushy whiskers and is universally respected, is apt to walk on tip-toe;
  \item I delight in cutting a man, whom I perfectly remember, but who will not give up his umbrella to me when it is raining. 
\end{enumerate}}
This problem can be rewritten using 52 variables into 50 nullities, which we take from \cite{GH1}.

{\tiny
\begin{tabular}{ll}
  $A="\text{London-friends of mine}"$&
  $B="\text{magistrates}"$\\
  $C="\text{men to whom I try to be civil}"$&
  $D="\text{men who cut me}"$\\
  $E="\text{men who have lent me money}"$&
  $F="\text{men whom I dare insult}"$\\
  $G="\text{men whom I delight to cut}"$&
  $H="\text{men with whom I like to walk}"$\\
  $I="\text{willing to shake hands with me when I am in rags}"$&
  $K="\text{novelists}"$\\
  $J="\text{men with whom the Colonel is willing to play billiards}"$&
  $L="\text{old}"$\\
  $M="\text{regarded by me with contemptuous pity}"$&
  $N="\text{remembered by me}"$\\
  $O="\text{yielding to me in argument}"$&
  $P="\text{rich}"$\\
  $Q="\text{rousing me to enthusiasm}"$&
  $R="\text{running across a street}"$\\
  $S="\text{sometimes losing his umbrella}"$&
  $T="\text{sticking to business}"$\\
  $U="\text{subjects of general conversation}"$&
  $V="\text{successful in life}"$\\
  $W="\text{tight-rope dancers}"$&
  $X="\text{universally respected}"$\\
  $Y="\text{wearing a coat}"$&
  $Z="\text{wearing a tall hat}"$\\
  $a="\text{apt to tell pointless anecdotes}"$&
  $b="\text{apt to walk on tip-toe}"$\\
  $c="\text{`at home` to me}"$&
  $d="\text{attractive}"$\\
  $e="\text{caring for appearances}"$&
  $f="\text{choosing his own wines}"$\\
  $g="\text{dull companions}"$&
  $h="\text{early risers}"$\\
  $i="\text{wearing kid gloves}"$&
  $j="\text{easily taken in}"$\\
  $k="\text{elephant-hunters}"$&
  $l="\text{farmers}"$\\
  $m="\text{fat}"$&
  $n="\text{good shots}"$\\
  $o="\text{willing to give up umbrella to me when it is raining}"$&
  $p="\text{good-tempered}"$\\
%
  $q="\text{having bushy whiskers}"$&
  $r="\text{having his wits about him}"$\\
  $s="\text{horsey}"$&
  $t="\text{humbugs}"$\\
  $u="\text{intimate friends of mine}"$&
  $v="\text{invited to dine with me}"$\\
  $w="\text{knowing what o'clock it is}"$&
  $x="\text{likely to be elected Mayor}"$\\
  $y="\text{likely to get a bad fall}"$&
  $z="\text{likely to marry an heiress}"$\\[1em]
\end{tabular}}

The nullities then are given by:
{\footnotesize
\begin{multicols}{4}   
\begin{enumerate}
  \item $HAZ'  =0$;
  \item $J'U'm'=0$;
  \item $wO'd  =0$;
  \item $WLQ   =0$;
  \item $KGD'  =0$;
  \item $iY't' =0$;
  \item $uAL   =0$;
  \item $UQs'  =0$;
  \item $rf'c' =0$;
  \item $na'p' =0$;
  \item $tKu'  =0$;
  \item $HpY'  =0$;
  \item $S'j'h'=0$;
  \item $mW'X' =0$;
  \item $M'V'R =0$;
  \item $T'q'x =0$;
  \item $kQ'l' =0$;
  \item $Aat'  =0$;
  \item $LvE'  =0$;
  \item $T'R'qy=0$;
  \item $h'Sz  =0$;
  \item $Lei'  =0$;
  \item $pt'D  =0$;
  \item $a'E'r'=0$;
  \item $fOv'  =0$;
  \item $C'V'q'=0$;
  \item $ls'   =0$;
  \item $Kg'Q' =0$;
  \item $hTXP' =0$;
  \item $ACx'  =0$;
  \item $pEe'I'=0$;
  \item $Js'l' =0$;
  \item $dPv'  =0$;
  \item $bMT'  =0$;
  \item $cv'B' =0$;
  \item $u'Io  =0$;
  \item $AsX'  =0$;
  \item $d'fj' =0$;
  \item $KF'n' =0$;
  \item $Lby'  =0$;
  \item $w'E'z'=0$;
  \item $Arj   =0$;
  \item $LN'I  =0$;
  \item $KT'V  =0$;
  \item $p'Fu' =0$;
  \item $BI'f' =0$;
  \item $gs'k' =0$;
  \item $ZiD'  =0$;
  \item $qXb'  =0$;
  \item $G'No' =0$.
\end{enumerate}
\end{multicols}
}

Keep in mind that Lewis Carroll's approach would get us with a tremendous amount of work a nullitity yielding in effect a lower bound, but not the "best", i.e. largest lower bound. This logical puzzle is indeed too big to be solved completely -- in our more restrictive sense -- by hand calculation, even the naive implementation of the 
solution routine presented here fails. 
However, using the SAT-structure and a recursive DPLL (Davis-Putnam-Logemann-Loveland)
SAT solver, see \cite{Davis62}, it can actually be done in less than two seconds on a standard laptop. 
The output (without the list of eliminands) is:

\begin{minipage}[t]{\textwidth}
\begin{lstlisting}
>> LewisCarrollLogic
Example 40: Friends
retinends: (A), (H), (K), (L), not (Y)
conclusion: 0 = (A) and (K) and (L) and not (Y) 
            0 = (A) and (H) and (K) and (L) 
\end{lstlisting}
\end{minipage} 
Thus, there is a unique proper, complete conclusion
$$\max\{AKLY^\prime,AKLH\}=0.$$
This could be verbalized as  "my old novelist friends from London all either wear coats or are the kind of men I do not like to walk with".

\end{document}